\documentclass{amsart}

\newtheorem{theorem}{Theorem}[section]
\newtheorem{lemma}[theorem]{Lemma}
\usepackage{color}
\theoremstyle{definition}
\newtheorem{definition}[theorem]{Definition}

\newtheorem{example}[theorem]{Example}

\usepackage{graphicx,subfigure}
\theoremstyle{remark}
\newtheorem{remark}[theorem]{Remark}

\numberwithin{equation}{section}

\usepackage{epstopdf}

\newcommand{\uu}{\textbf{{u}}}

\newcommand{\f}{\textbf{{f}}}

\begin{document}

\title[Well-Conditioned Galerkin Method for Two-Sided Fractional Equation]{Well-Conditioned Galerkin Spectral Method for Two-Sided Fractional Diffusion Equation with Drift}

\author{Lijing Zhao}
\address{School of Natural and Applied Sciences, Northwestern Polytechnical University,
Xi'an 710129, P.R. China}
\email{zhaolj17@nwpu.edu.cn}
\thanks{The first author is supported by NSFC Grant 11801448, by the Natural Science Basic Research Plan in Shaanxi Province of China under Grant 2018JQ1022.}

\author{Xudong Wang}
\address{School of Mathematics and Statistics, Gansu Key Laboratory of Applied Mathematics and Complex Systems, Lanzhou University, Lanzhou 730000, P.R. China}
\email{xdwang14@lzu.edu.cn}
\thanks{The second author is supported in part by the Fundamental Research Funds for the Central Universities under Grants No. lzujbky-2019-it17.}

%





\begin{abstract}
In this paper, we focus on designing a well-conditioned Glarkin spectral methods for solving a two-sided fractional diffusion equations with drift, in which the fractional operators are defined neither in Riemann-Liouville nor Caputo sense, and its physical meaning is clear.
Based on the image spaces of Riemann-Liouville fractional integral operators on $L_p([a,b])$ space discussed in our previous work, after a step by step deduction, three kinds of Galerkin spectral formulations are proposed,  the final obtained corresponding scheme of which shows to be well-conditioned---the condition number of the stiff matrix can be reduced from $O(N^{2\alpha})$ to $O(N^{\alpha})$, where $N$ is the degree of the polynomials used in the approximation. Another point is that the obtained schemes can also be applied successfully to approximate fractional Laplacian with generalized homogeneous boundary conditions, whose fractional order $\alpha\in(0,2)$, not only having to be limited to $\alpha\in(1,2)$.
Several numerical experiments demonstrate the effectiveness of the derived schemes. Besides, based on the numerical results, we can observe the behavior of mean first exit time, an interesting quantity that can provide us with a further understanding about the mechanism of abnormal diffusion.
\end{abstract}

\maketitle

%

\section{Introduction}\label{section:1}
In this paper, we target on investigating a well-condtioned Galerkin spectral methods for the following two-sided fractional diffusion equation with drift
\begin{equation}\label{main_equation}
\left\{ \begin{array}{l}
-\left(\bar{p}\cdot\textbf{D}^{\alpha}u(x)+\bar{q}\cdot\textbf{D}^{\alpha \ast}u(x)\right)
+\bar{d}\cdot{}Du(x)=h(x),~~x\in(a,b),\\
u(a)=u(b)=0,
\end{array} \right.
\end{equation}
where $1<\alpha< 2$, $0\leq \bar{p},\bar{q}\leq 1$ satisfying $\bar{p}+\bar{q}=1$, and $\textbf{D}^{\alpha}$, $\textbf{D}^{\alpha \ast}$ are neither the Riemann-Liouville operators nor the Caputo ones \cite{Podlubny:99}; rather, in general, for $n-1\leq\gamma<n$,
\begin{equation}\label{left_def}
\textbf{D}^{\gamma}u(x):=D{}_{a}I_{x}^{n-\gamma}D^{n-1} u(x),
\end{equation}
\begin{equation}\label{right_def}
\textbf{D}^{\gamma \ast}u(x):=(-1)^{n}D{}_{x}I_{b}^{n-\gamma}D^{n-1} u(x),
\end{equation}
with $_{a}I_{x}^{\beta}$ and $_{x}I_b^{\beta}$, $0<\beta<1$, denote separately the left Riemann-Liouville fractional integral
$$_{a}I_{x}^{\beta}u(x)=\frac{1}{\Gamma(\beta)}\int_{a}^x(x-s)^{\beta-1}u(s)ds,$$
and the right Riemann-Liouville fractional integral operator
$$_{x}I_b^{\beta}u(x)=\frac{1}{\Gamma(\beta)}\int_{x}^b(s-x)^{\beta-1}u(s)ds.$$
Here $\Gamma(\cdot)$ presents the Euler gamma function.

The fractional Dirichlet problem and variants thereof appear in many applications, in particular in physical settings where anomalous dynamics occur and where the spread of mass grows faster than linearly in time. Examples include turbulent fluids, contaminant transport in fractured rocks, chaotic dynamics and disordered quantum ensembles; see \cite{Klafter:11,Klages:08,Shlesinger:95}. The authors in \cite{Ervin:16} believe that problem (\ref{main_equation}), which can be interpreted as the steady-state equation for a time dependent advection and anomalous diffusion problem, is a more physical model than the corresponding Riemann-Liouville or Caputo fractional equation. During the derivation of Eq. (\ref{main_equation}), the authors in \cite{Ervin:16} point out that besides obeying the conservation of energy principle, the physical interpretation of the flux at a given cross section $x$, is that ``there is a nonlocal effect from a flux originating at a cross section $s$, proportional to $1/($distance for that point$)^{\alpha-1}$". In other words, the contribution to the flux at cross section $x$, from points to its left and right is given by
$$-k\int_{a}^x (x-s)^{1-\alpha}\frac{\partial u(s,t)}{\partial s}\,ds,$$
and
$$-k\int_{x}^b (s-x)^{1-\alpha}\frac{\partial u(s,t)}{\partial s}\,ds,$$
respectively, where $k$ is a dispersion coefficient.
In this way, when considering the case of Dirichlet boundary conditions $u(a)=u(b)=1$, after by changing the unknown $v=u-1$, the simulation of the model equation would require the same energy source as for the case $u(a)=u(b)=0$, which physically makes sense.

Besides can be viewed as the steady-state equation for a time dependent advection and anomalous diffusion problem, we shall see that when $h(x)=-1$, under the framework of the image spaces of Riemann-Liouville fractional integral operators on $L_p([a,b])$ space \cite{Zhao:16}, Problem (\ref{main_equation}) itself can also be used to describe the mean first exit time of a stochastic process never leaving a fixed region in the state space \cite{Deng:17}---an interesting deterministic quantity that can provide us with a further understanding about the mechanism of the anomalous diffusion.

Another topic we want to note is that from a mathematical view, under suitable assumptions on $u$, fractional Laplacian operator \cite{Yang:10}
\begin{equation}\label{laplacian_operator}
(-\triangle)^s u(x)=\frac{2^{2s-1}\Gamma(s+\frac{1}{2})}{\pi^{\frac{1}{2}}\Gamma(1-s)}
\int_{\mathbb{R}}\frac{u(x+y)+u(x-y)-2u(x)}{|y|^{1+2s}}\,dy.
\end{equation}
is equivalent to
\begin{equation*}
\begin{array}{rll}
(-\triangle)^{\frac{\alpha}{2}} u(x)
&=&\frac{1}{2\cos(\frac{\pi \alpha}{2})}
\left({}^{RL}_{-\infty}D_{x}^{\alpha}+{}^{RL}_{x}D_{\infty}^{\alpha}\right)u(x)\\[5pt]
£º&=&\frac{1}{2\cos(\frac{\pi \alpha}{2})}
\left(D^2{}_{-\infty}I_{x}^{\alpha}+D^2{}_{x}I_{\infty}^{\alpha}\right)u(x),
\end{array}
\end{equation*}
where \textbf{$\alpha\in(0,2)$}.
Actually, the proof in \cite{Yang:10} also ensures that
\begin{equation*}
(-\triangle)^{\frac{\alpha}{2}} u(x)
=\frac{1}{2\cos(\frac{\pi \alpha}{2})}
\left(D{}_{-\infty}I_{x}^{\alpha}D+D{}_{x}I_{\infty}^{\alpha}D\right)u(x).
\end{equation*}
Therefore, mathematically, the following one-dimensional Poisson problem with generalized Dirichlet boundary condition:
\begin{equation}\label{equation_laplacian}
\left\{ \begin{array}{rlll}
(-\triangle)^{\alpha/2}u(x)&=&h(x),&~x\in\Omega=(a,b),\\[5pt]
u(x)&=&0,&~x\in \mathbb{R}\backslash  \Omega,
\end{array} \right.
\end{equation}
can be changed as
\begin{equation}\label{laplace_to_fractional}
\left\{ \begin{array}{lll}
-\frac{1}{2}\left(\textbf{D}^{\alpha}u(x)+\textbf{D}^{\alpha \ast}u(x)\right)=-\cos(\frac{\pi \alpha}{2}) h(x),~~x\in(a,b),\\
u(a)=u(b)=0,
\end{array}\right.
\end{equation}
which is a special case of (\ref{main_equation}), where \textbf{$\alpha\in(0,2)$} \textbf{(not only limited to (1,2) as in (\ref{main_equation}))}.

From the view of stochastic processes, the physical meaning of the fractional Laplacian defined in above way with Dirichlet boundary conditions is the negative infinitesimal generator of stopped subordinated Brownian motion (i.e., stopped $\alpha$-stable L\'{e}vy motion), which represents particles that are stopped upon exiting the domain via a jump over the boundary \cite{Deng:17,Deng:17-2,Meerschaert:11}.
Here, we do not concern the detailed conditions under which (\ref{equation_laplacian}) and (\ref{laplace_to_fractional}) are equivalent. Instead, we mainly focus on the spectral methods that are effective for them, and leave the theoretical part in our future work.


Comparing with the classical differential equations, one of the big challenges we have to face is the expensiveness of its computation cost besides its complexity, since fractional operators
are pseudodifferential operators which are non-local. Finite difference methods and finite elements methods are not easy to apply when solving especially a two-sided fractional problems, because the information on the whole domain is needed which results in a huge computational cost. In this case, spectral method, as a global method, appears to be a natural choice. There are existing spectral work, used to solve one-sided or two-sided fractional differential equations with Riemann-Liouville or variable order fractional operators \cite{Chen:16,Li:12,Mao:18,Wang:15,Zayernouri:15}.
Early spectral collocation methods for fractional problems using classical
interpolation basis functions with
Legendre-Gauss-Lobatto or Chebyshev-Gauss-Lobatto collocated points are proposed
in \cite{Li:09} and \cite{Tian:14}. Eigenfunctions of a fractional Sturm-Liouville operator are derived in \cite{Zayernouri:13}. Spectral approximation results in weighted Sobolev spaces involving
fractional derivatives are derived in \cite{Chen:16}, including also rigorous convergence analysis. The authors in \cite{Jiao:16} introduce fractional Birkhoff interpolation basis
functions into collocation methods to reduce the condition numbers when solving the one-sided Caputo fractional equations.

As for the problem (\ref{main_equation}), a variational formulation is studied in \cite{Ervin:06}, together with a finite element error analysis. The regularity of (\ref{main_equation}) is studied, also a finite elements method and a spectral type approximation method are proposed in it.

As far as we know, there is little literature to discuss the weak formulation of the two-sided fractional diffusion problems with drift, in which the fractional operators are physically well-defined. Also, there has been no relevant work to talk about the corresponding well-conditioned scheme.

This paper mainly proposes three kinds of Galerkin spectral schemes for solving Eq. (\ref{main_equation}). These three Galerkin spectral schemes are based on different weak variational formulations and have different regularity requirements, all of which shows to be effective to this kind of two-sided fractional diffusion equation with drift, even when the solution has a low regularity. In special, based on the former two formulations, the third one, named as mixed Galerkin spectral formulation, is designed by splitting the Eq. (\ref{main_equation}) into three subequations. In this way, the trial and test functions are more flexible to choose, so that the coefficient matrices can be expressed in a simpler way. Besides, compared with the condition numbers of the stiff matrices in the other two schemes, the condition number in mixed Galerkin spectral scheme can be reduced from $O(N^{2\alpha})$ to about $O(N^{\alpha})$, where $N$ is the degree of the polynomials used in the approximation.

The rest of this paper is organized as follows. Section \ref{section:2} reviews some important definitions and results about the image spaces of Riemann-Liouville fractional integrals on Sobolev space $W^{m,p}(\Omega)$ space, which are the framework of the weak formulation in this paper.
Three different weak formulations and Galerkin spectral methods are presented step by step in Sections \ref{section:3}, where the differences among them are discussed. Section \ref{section:4} provides
the numerical results for solving problems (\ref{main_equation}) and (\ref{equation_laplacian}), in which one can observed that the condition numbers are substantially decreased in the mixed Galerkin spectral method. Finally, the main results are summarized in Section \ref{section:5}.

\section{Preliminaries }\label{section:2}

In this section, we outline the definition and some results about the image spaces of $\alpha$-order Riemann-Liouville fractional integral operators on $L_p(\Omega)$ or $W^{m,p}(\Omega)$, which is called ``spaces of fractional integrals" for short \cite{Zhao:16}, where $W^{m,p}(\Omega)$ is a given classical Sobolev space and $\Omega=[a,b]$.

As we all know that the concept of fractional calculus is almost as old as their more familiar integer order counterparts, and many mathematical results about fractional operators are also discussed in the early days \cite{Miller:93,Oldham:74,Podlubny:99,Samko:93}. Until recently, fractional derivatives have been widely and successfully explored as a tool for developing
more sophisticated mathematical models. Here, we borrow (not simply copy but sometimes have to
flip through pages) the space, which we call as the image space of Riemann-Liouville fractional integral operators on $L_p([a,b])$ space, introduced in \cite{Samko:93}, and some results given in \cite{Samko:93} and \cite{Zhao:16}, to begin our discussion.
The reason we choose this kind of space, not only because it comes from a ``non state of the art" references, but also because the key difficulty of the fractional operators that are widely used, such as Riemann-Liouville derivative or Caputo deriavative, are actually come from the pseudo-differential or Riemann-Liouville fractional integral operator in them. Since the space of fractional integrals of $L_p$ functions can catch this characteristic very well, it is a natural way to begin our discussion from it.

Denote $L_p(\Omega)~ (1\le p<\infty)$ as $L_p$ space on $\Omega=[a,b]$. The set of $\alpha$-th order left and right Riemann-Liouville fractional integrals of $L_p(\Omega)$ functions, $1\le p<\infty$, are firstly given in Definition 2.3 of \cite{Samko:93}. We rearrange them as follows:
\begin{definition}\label{defi_frac_int_space}
\begin{equation}\label{defi_left}
I^{\alpha}[L_p(\Omega)]:=\left\{f:f(x)={}_{a}I_{x}^\alpha \varphi(x), \varphi(x)\in L_p(\Omega), x\in\Omega\right\},~\alpha>0,
\end{equation}
and
\begin{equation}\label{defi_right}
I^{\alpha\ast}[L_p(\Omega)]:=\left\{f:f(x)={}_{x}I_{b}^\alpha \varphi(x), \varphi(x)\in L_p(\Omega), x\in\Omega\right\},
~\alpha>0.
\end{equation}
\end{definition}

Now we only list some results about $I^{\alpha}[L_p(\Omega)]$; similar results can be derived for $I^{\alpha\ast}[L_p(\Omega)]$.

In \cite{Zhao:16}, Corollary 2.10 shows that actually if $u(x)\in I^{\alpha}[L_p(\Omega)]$, then ${}^{RL}_{a}D_{x}^{\alpha}u(x):=D^n{}_{a}I_{x}^{n-\alpha}u(x)=\textbf{D}^{\alpha}u(x)$.
Therefore, the following lemmas hold \cite{Zhao:16}.
\begin{lemma}\label{lemma1}
If $u(x)\in I^{\alpha}[L_1(\Omega)]$, $n-1\leq \alpha<n$, then
\begin{equation}\label{equation1_lemma1}
{}_{a}I_{x}^{\alpha}\textbf{D}^{\alpha}u(x)
=u(x).
\end{equation}
\end{lemma}

\begin{lemma}\label{lemma4}
Let $n-1\leq \alpha<n$. If $u(x)\in I^{\alpha}[L_2(\Omega)]$,
$v(x)\in I^{\alpha\ast}[L_2(\Omega)]$, then
\begin{equation}\label{equation:lemma4}
\left(\textbf{D}^{\alpha}u(x),v(x)\right)=\left(u(x),\textbf{D}^{\alpha\ast}v(x)\right).
\end{equation}
\end{lemma}

\begin{lemma}\label{lemma6}
Let $\alpha_1>0$, $\alpha_2>0$, $\alpha_1+\alpha_2=\alpha$. If
$u(x)\in I^{\alpha}[L_p(\Omega)]$, then
\begin{equation}\label{equation_1_lemma6}
\textbf{D}^{\alpha}u(x)=\textbf{D}^{\alpha_1}\textbf{D}^{\alpha_2}u(x),
\end{equation}
and
\begin{equation}\label{equation_2_lemma6}
\textbf{D}^{\alpha_2}u(x)\in I^{\alpha_1}[L_p(\Omega)].
\end{equation}
\end{lemma}

If $u(x)\in I^{\alpha}[L_1(\Omega)]$, then there exists a unique $\varphi(x)\in L_1(\Omega)$ \cite{Samko:93,Zhao:16}, such that $u(x)={}_{a}I_{x}^{\alpha} v(x)$. Using Lemma \ref{lemma1}, we have
\begin{equation}\label{weak_solution}
\begin{array}{lll}
&&\int_a^b u(x)\cdot \textbf{D}^{\alpha\ast}\phi(x)\,dx\vspace{3pt}\\
&=&\int_a^b {}_{a}I_{x}^{\alpha} \varphi(x)\cdot \textbf{D}^{\alpha\ast}\phi(x)\,dx\vspace{3pt}\\
&=&\int_a^b \varphi(x)\cdot {}_{x}I_{b}^{\alpha}{}\textbf{D}^{\alpha\ast}\phi(x)\,dx\vspace{3pt}\\
&=&\int_a^b \varphi(x)\phi(x)\,dx \qquad \qquad \qquad\forall \phi(x)\in C^{\infty}_{c}(\Omega),
\end{array}
\end{equation}
where the integrals make sense because of the H\"{o}lder inequality $\|fg\|_{L_1}\leq \|f\|_{L_1}\cdot\|g\|_{L_{\infty}}$.

Because $I^{\alpha}[L_p(\Omega)]\subseteq I^{\alpha}[L_1(\Omega)]$, $p\geq 1$, so, Eqs. (\ref{equation1_lemma1})-(\ref{weak_solution}) still hold for $I^{\alpha}[L_p(\Omega)]$, $p\geq 1$. Therefore, we can say that $I^{\alpha}[L_p(\Omega)]$ is a Sobolev space.

Since for $p\geq 1$,
$I^{\alpha}[L_p(\Omega)]\hookrightarrow L_p(\Omega)$ (Theorem 2.6 in \cite{Samko:93}), i.e., (similar to Poincar\'{e} inequality \cite{Ervin:06})
\begin{equation*}
\|{}_{a}I_{x}^\alpha\varphi(x)\|_{p}\le \frac{(b-a)^\alpha}{\Gamma(\alpha+1)}\|\varphi(x)\|_{p}
\qquad ~~~\forall \varphi(x)\in L_p(\Omega).
\end{equation*}
We can introduce the norm in $I^{\alpha}[L_p(\Omega)]$ by
\begin{equation}\label{norm}
\|u(x)\|_{I^{\alpha}[L_p(\Omega)]}:=\|{}_{a}\textbf{D}_{x}^{\alpha}u(x)\|_{p}.
\end{equation}


\begin{remark}\label{remark1}
In the later sections, we can see that actually, for $\alpha>0$, $\delta>-1$, $\gamma\in \mathbb{R}$, functions $(1+x)^{\delta+\alpha} J_{n}^{\gamma-\alpha,\delta+\alpha}(x)$ and
$(1-x)^{\delta+\alpha} J_{n}^{\delta+\alpha,\gamma-\alpha}(x)$ belong to $I^{\alpha}[L_1(-1,1)]$ and $I^{\alpha\ast}[L_p(-1,1)]$, respectively, where $\{J_{n}^{\sigma,\eta}(x)\}_{n=0}$ denote the Jacobi polynomials, which are defined by Rodrigues' formula
\begin{equation*}
(1-x)^\sigma(1+x)^\eta J_{n}^{\sigma,\eta}(x)
=\frac{(-1)^n}{2^n n!}\frac{d^n}{dx^n}\left[(1-x)^{n+\sigma}(1+x)^{n+\eta}\right],
\end{equation*}
and they are orthogonal on $[-1,1]$ with respect to
$(1-x)^\sigma(1+x)^\eta$ when $\sigma>-1$, $\eta>-1$ \cite{Szego:75}.
\end{remark}

Next, the Sobolev space with higher regularity can be defined \cite{Zhao:16}:
\begin{definition}\label{defi_higher_space}
The image space of $\alpha$-th order left Riemann-Liouville fractional integrals on $W^{m,p}(\Omega)$ is defined as
\begin{equation}\label{defi_left_higher}
I^\alpha\left[W^{m,p}(\Omega)\right]:=
\left\{f:f(x)={}_{a}I_{x}^\alpha \varphi(x), \varphi(x)\in W^{m,p}(\Omega), x\in\Omega\right\},
\end{equation}
and with norm
\begin{equation*}
\|f(x)\|_{I^{\alpha}[W^{m,p}(\Omega)]}:=
\|\textbf{D}_{x}^{\alpha}f(x)\|_{W^{m,p}(\Omega)},
\end{equation*}
where $W^{m,p}(\Omega)$ is a given classical integer Sobolev space.
\end{definition}

The relationships between the image spaces $I^\alpha\left[W^{m,p}(\Omega)\right]$ and $I^{\alpha\ast}\left[W^{m,p}(\Omega)\right]$ are briefly listed in the following lemmas, in which besides the case $m=0$, the most interested cases is when $p=2$ and $W^{m,2}(\Omega)=H^{m}(\Omega)$.

\begin{lemma}\cite{Samko:93}\label{lemma2}
When $0<\alpha<1/p$, and $1<p<\infty$, then
\begin{equation}\label{left_right_Lp_1}
H^{\alpha,p}(\Omega)=\hat{I}^{\alpha}[L_p(\Omega)]:=
I^{\alpha}[L_p(\Omega)]=I^{\alpha\ast}[L_p(\Omega)].
\end{equation}

When $1/p<\alpha<1/p+1$, then
\begin{equation}\label{left_right_Lp_2}
H_0^{\alpha,p}(\Omega)=I^{\alpha}[L_p(\Omega)]\cap I^{\alpha\ast}[L_p(\Omega)],
\end{equation}
where
\begin{equation*}
H_0^{\alpha,p}(\Omega)=\left\{f:f(x)\in H^{\alpha,p}(\Omega), \textrm{~and~} f(a)=f(b)=0\right\},
\end{equation*}
\begin{equation*}
H^{\alpha,p}(\Omega)=\left\{f:\exists ~g(x)\in H^{\alpha,p}(\mathbb{R}),
~\textrm{s.t.}~ g(x)\big|_{\Omega}=f(x)\right\},
\end{equation*}
\begin{equation*}
H^{\alpha,p}(\mathbb{R})=\{f(x)\in L_p(\mathbb{R}):\mathcal{F}^{-1}[(1+|\xi|^2)^{\frac{\alpha}{2}}
\mathcal{F}[f] ]  \in L_p(\mathbb{R})\}.
\end{equation*}
\end{lemma}

\begin{lemma}\cite{Zhao:16}\label{lemma3}
If ${0\leq\alpha<\frac{1}{2}}$, then
\begin{eqnarray}\label{left_right_Wmp_1}
&&I^{\alpha}[H^{m}(\Omega)]\cap I^{\alpha\ast}[H^{m}(\Omega)]
\nonumber\\
&=&\bigg\{
f:f(x)\in W^{{m},q}(\Omega),
f(x)=o((x-a)^{m+\alpha-\frac{1}{2}}),~\textrm{as}~x\rightarrow a,
\nonumber\\
&&\quad
f(x)=o((b-x)^{m+\alpha-\frac{1}{2}}),~\textrm{as}~x\rightarrow b\bigg\},
~~{q=\frac{2}{1-2\alpha}}.
\end{eqnarray}

If ${\frac{1}{2}<\alpha<1}$, then
\begin{eqnarray}\label{left_right_Wmp_2}
&&I^{\alpha}[H^{m}(\Omega)]\cap I^{\alpha\ast}[H^{m}(\Omega)]
\nonumber\\
&=&\bigg\{
f:f(x)\in W^{{m+1},q}(\Omega),
f(x)=o((x-a)^{m+\alpha-\frac{1}{2}}),~\textrm{as}~x\rightarrow a,
\nonumber\\
&&\quad
f(x)=o((b-x)^{m+\alpha-\frac{1}{2}}),~\textrm{as}~x\rightarrow b\bigg\},
~~{q=\frac{2}{3-2\alpha}}.
\end{eqnarray}
\end{lemma}

Denote $P_{N}(\Omega)$ as the polynomials spaces of degree less than or equal to $N$ on $\Omega$. Then $I^{\alpha}\left[P_{N}(\Omega)\right]:=\left\{f:f(x)={}_{a}I_{x}^\alpha \varphi(x), \varphi(x)\in P_{N}(\Omega), x\in[\Omega]\right\}$ is a subspace of $ I^{\alpha}[L_2(\Omega)]$.

Denote $\Pi_N$ as the orthogonal projection operator from  $L_2(\Omega)$ onto $P_{N}(\Omega)$. Then the following
approximation property holds:
\begin{lemma}\cite{Zhao:16}\label{lemma5}
If $\alpha\in(0,\frac{1}{2})\cup (\frac{1}{2},1)$, and $u(x)\in I_{a+}^{\alpha}[H^m(\Omega)]$, then there exists a constant $C=C(\alpha,\Omega,m)$, such that
\begin{equation}\label{equation:lemma5}
\|u-Q_N^{\alpha} u\|_{L_2(\Omega)}\leq C N^{-m}\|\,_{a}\textbf{D}_{x}^{\alpha} u\|_{H^m(\Omega)},
\end{equation}
where $Q_N^\alpha u(x):=\,_{a}I_{x}^{\alpha}\left(\Pi_N \,_{a}\textbf{D}_{x}^{\alpha}u\right)(x)$.
\end{lemma}
\section{Variational formulations and spectral methods}\label{section:3}

We use the spaces of fractional integrals introduced above to design Galerkin spectral methods for solving problem (\ref{main_equation}). Without loss of generality, we now restrict our attention to the interval $\Omega=[-1,1]$.

\subsection{Variational formulations}

In order to derive a variational form of (\ref{main_equation}),
we firstly assume for the moment that $u(x)$ is a sufficiently smooth solution.
By multiplying an arbitrary $v(x)\in C_{c}^{\infty}(\Omega)$, it can be obtained that
\begin{equation}\label{equation_integrate}
\int_{-1}^{1}-\left(\bar{p}\cdot\textbf{D}^{\alpha}u(x)+\bar{q}\cdot\textbf{D}^{\alpha \ast}u(x)\right)\cdot v(x)\,dx
+\bar{d}{}Du(x)\cdot v(x)\,dx
=\int_{-1}^{1}h(x)v(x)\,dx.
\end{equation}

\subsubsection{Variational formulation \uppercase\expandafter{\romannumeral1}}\label{subsbusection:3.1.1}
Taking integration by parts for the left hand of (\ref{equation_integrate}), and noting that $\textbf{D}u(x)=\textbf{D}^{\frac{\alpha}{2}}\textbf{D}^{\frac{\alpha}{2}}u(x)$ for smooth $u$ with $u(-1)=0$ by the definition of $\textbf{D}^{\gamma}$ in(\ref{left_def}), we can obtain
\begin{equation}\label{equation_weak_integration1}
\begin{array}{lll}
&&-\bar{p}\int_{-1}^{1}\textbf{D}^{\frac{\alpha}{2}}u(x)\cdot\textbf{D}^{\frac{\alpha}{2}\ast}v(x)\,dx
-\bar{q}\int_{-1}^{1}\textbf{D}^{\frac{\alpha}{2}\ast}u(x)\cdot\textbf{D}^{\frac{\alpha}{2}}v(x)\,dx\\
&&+\bar{d}\int_{-1}^{1}\textbf{D}^{\frac{1}{2}}u(x)\cdot\textbf{D}^{\frac{1}{2}\ast}v(x)\,dx
=\int_{-1}^{1}h(x)v(x)\,dx.
\end{array}
\end{equation}

Denote
$$\Phi_{1}^{\frac{\alpha}{2}}(\Omega):=I^{\frac{\alpha}{2}}[L_2(\Omega)]\cap I^{\frac{\alpha}{2}\ast}[L_2(\Omega)].$$
Now we define the associated bilinear form
$B_1:\Phi_{1}^{\frac{\alpha}{2}}(\Omega)\times \Phi_{1}^{\frac{\alpha}{2}}(\Omega)\rightarrow \mathbb{R}$ for (\ref{main_equation}) as
\begin{equation}\label{bilinear_form_1}
B_1(u,v):=-\bar{p}\left(\textbf{D}^{\frac{\alpha}{2}}u,\textbf{D}^{\frac{\alpha}{2}\ast}v\right)
-\bar{q}\left(\textbf{D}^{\frac{\alpha}{2}\ast}u,\textbf{D}^{\frac{\alpha}{2}}v\right)
+\bar{d}\cdot\left(\textbf{D}^{\frac{1}{2}}u,\textbf{D}^{\frac{1}{2}\ast}v\right).
\end{equation}

For a given function $h(x)$, which belongs to the dual space of $W_{0}^{1,p_1}(\Omega)$ \cite{Brezis:10}, and be denoted as $W^{-1,q_1}(\Omega)$, where $p_1=\frac{2}{3-\alpha}$, $q_1=\frac{2}{\alpha-1}$, we define the associated linear functional $F_1:\Phi_{1}^{\frac{\alpha}{2}}(\Omega)\rightarrow \mathbb{R}$ as
\begin{equation}\label{functional_form_1}
F_1(v):=\langle h,v\rangle,
\end{equation}
where $\langle \cdot,\cdot \rangle$ is the duality pair of $W^{-1,q_1}(\Omega)$ and $W_{0}^{1,p_1}(\Omega)$.

By Lemma \ref{lemma6} and formula (\ref{left_right_Wmp_2}) in Lemma \ref{lemma3}, we can check that both (\ref{bilinear_form_1}) and (\ref{functional_form_1}) make sense.

Thus, the corresponding variational formulation of (\ref{main_equation}) can be defined as follows.

\begin{definition}[\textrm Variational Formulation \uppercase\expandafter{\romannumeral1}]\label{weak_form_1}
A function $u(x)\in \Phi_{1}^{\frac{\alpha}{2}}(\Omega)$ is a variational solution of problem
(\ref{main_equation}) provided
\begin{equation}\label{equation:weak_form_1}
B_1(u,v)=F_1(v)\quad\forall v(x)\in \Phi_{1}^{\frac{\alpha}{2}}(\Omega).
\end{equation}
\end{definition}

Denote
$$\Phi_{1,N}^{\frac{\alpha}{2}}(\Omega)=I^{\frac{\alpha}{2}}[P_N(\Omega)]\cap I^{\frac{\alpha}{2}\ast}[P_N(\Omega)].$$
Then the Galerkin approximation of (\ref{equation:weak_form_1}) is: find $u_{1,N}(x)\in \Phi_{1,N}^{\frac{\alpha}{2}}(\Omega)$, such that
\begin{equation}\label{equation:Galerkin_form_1}
B_1(u_{1,N},v_{1,N})=F_1(v_{1,N})\quad\forall v_{1,N}(x)\in \Phi_{1,N}^{\frac{\alpha}{2}}(\Omega).
\end{equation}

\subsubsection{Variational formulation \uppercase\expandafter{\romannumeral2}}\label{subsbusection:3.1.2}

Actually, for smooth solution $u$ with $u(-1)=u(1)=0$, and an arbitrary given $v(x)\in C_{c}^{\infty}(\Omega)$, instead of Eq. (\ref{equation_weak_integration1}), we can get another formula by taking integration by part for the left side of Eq. (\ref{equation_integrate}), as follows:
\begin{equation}\label{equation_weak_integration2}
\begin{array}{lll}
&&-\bar{p}\int_{-1}^{1}\textbf{D}^{\frac{\alpha-1}{2}}u(x)\cdot\textbf{D}^{\frac{\alpha+1}{2}\ast}v(x)\,dx
-\bar{q}\int_{-1}^{1}\textbf{D}^{\frac{\alpha-1}{2}\ast}u(x)\cdot\textbf{D}^{\frac{\alpha+1}{2}}v(x)\,dx\\[5pt]
&&~~~~~~~~~~~~-\bar{d}\int_{-1}^{1}u(x)\cdot D v(x)\,dx
=\int_{-1}^{1}h(x)v(x)\,dx.
\end{array}
\end{equation}

Denote $$\Phi_{2}^{\frac{\alpha-1}{2}}(\Omega):=\left\{f:f\in\hat{I}^{\frac{\alpha-1}{2}}[L_2(\Omega)], \textrm{~and~}f(-1)=f(1)=0\right\}.$$
We now define another type of bilinear form
$B_2:\Phi_{2}^{\frac{\alpha-1}{2}}(\Omega)\times \Phi_{1}^{\frac{\alpha+1}{2}}(\Omega)\rightarrow \mathbb{R}$ for (\ref{main_equation}) as
\begin{equation}\label{bilinear_form_2}
B_2(u,v):=-\bar{p}\left(\textbf{D}^{\frac{\alpha-1}{2}}u,\textbf{D}^{\frac{\alpha+1}{2}\ast}v\right)
-\bar{q}\left(\textbf{D}^{\frac{\alpha-1}{2}\ast}u,\textbf{D}^{\frac{\alpha+1}{2}}v\right)
-\bar{d}\cdot\left(u,Dv\right).
\end{equation}

For a given source term $h(x)$, which belongs to the dual space of $W_{0}^{1,p_2}(\Omega)$ \cite{Brezis:10}, and be denoted as $W^{-1,q_2}(\Omega)$, where $p_2=\frac{2}{2-\alpha}$, $q_2=\frac{2}{\alpha}$, we define the associated linear functional $F_2:\Phi_{1}^{\frac{\alpha+1}{2}}(\Omega)\rightarrow \mathbb{R}$ as
\begin{equation}\label{functional_form_2}
F_2(v):=\langle h,v\rangle,
\end{equation}
where $\langle \cdot,\cdot \rangle$ is the duality pair of $W^{-1,q_2}(\Omega)$ and $W_{0}^{1,p_2}(\Omega)$.

By Lemma \ref{lemma6}, formula (\ref{left_right_Lp_2}) in Lemma \ref{lemma2}, and formula (\ref{left_right_Wmp_1}) in Lemma \ref{lemma3}, we can check that both (\ref{bilinear_form_2}) and (\ref{functional_form_2}) make sense.

Thus, the corresponding variational formulation of (\ref{main_equation}) can be defined as follows.

\begin{definition}[\textrm Variational Formulation \uppercase\expandafter{\romannumeral2}]\label{weak_form_2}
A function $u(x)\in \Phi_{2}^{\frac{\alpha-1}{2}}(\Omega)$ is a variational solution of problem
(\ref{main_equation}) provided
\begin{equation}\label{equation:weak_form_2}
B_2(u,v)=F_2(v)\quad\forall v(x)\in \Phi_{1}^{\frac{\alpha+1}{2}}(\Omega).
\end{equation}
\end{definition}

\begin{remark}\label{remark2}
It is not difficulty to see that the weak solution as well as the linear functional in (\ref{equation:weak_form_2}) lie in weaker spaces than the weak solution and the linear functional of (\ref{equation:weak_form_1}) do; the classical solution can be recovered from both (\ref{equation:weak_form_1}) and (\ref{equation:weak_form_2}) if $u$ is smooth enough.
\end{remark}

Denote
$$\Phi_{2,N}^{\frac{\alpha-1}{2}}(\Omega):=\hat{I}^{\frac{\alpha-1}{2}}[P_N(\Omega)].$$
We can see that if $f(x)\in \Phi_{2,N}^{\frac{\alpha-1}{2}}(\Omega)$, then $f(\pm1)=0$.

The Galerkin approximation of (\ref{equation:weak_form_2}) is: find $u_{2,N}(x)\in \Phi_{2,N}^{\frac{\alpha-1}{2}}(\Omega)$, such that
\begin{equation}\label{equation:Petrov_Galerkin_form_2}
B_2(u_{2,N},v_{2,N})=F_2(v_{2,N})\quad\forall v_{2,N}(x)\in \Phi_{1,N}^{\frac{\alpha+1}{2}}(\Omega).
\end{equation}

\subsubsection{Variational formulation \uppercase\expandafter{\romannumeral3}}\label{subsbusection:3.1.3}

Since by Lemma \ref{lemma2}, when $\gamma>\frac{1}{2}$, $I^{\gamma}[P_N(\Omega)]\neq I^{\gamma\ast}[P_N(\Omega)]$, it is not simple to manipulate $I^{\gamma}[P_N(\Omega)]\cap I^{\gamma\ast}[P_N(\Omega)]$ during the numerical realization. One way to get rid of using it during the computation, is based on the following splitting formula, which is equivalent to problem (\ref{main_equation}):
\begin{equation}\label{split_equation}
\begin{split}
    \left\{ \begin{array}{ll}
        l(x)=\bar{p}~{}_{-1}I_x^{2-\alpha}Du(x), \\
        r(x)=\bar{q}~{}_xI_{1}^{2-\alpha}Du(x),  \\
        -D[l(x)+r(x)]+\bar{d}\cdot Du(x)=h(x),\\
        u(-1)=u(1)=0.
    \end{array}
  \right.
\end{split}
\end{equation}
Similarly to the above discussions, by assuming for the moment that $u(x)$ is a sufficiently smooth solution,
then multiplying the first three equalities of (\ref{split_equation}) separately by $\psi_1\in C_{c}^{\infty}(\Omega)$, $\psi_2\in C_{c}^{\infty}(\Omega)$, $\psi_3\in C_{c}^{\infty}(\Omega)$, and taking integration by parts, we can get
\begin{equation}\label{equation_weak_integration3}
\begin{split}
    \left\{ \begin{array}{ll}
        \int_{-1}^1 l(x)\psi_1(x)\,dx
        =\bar{p}~\int_{-1}^1\textbf{D}^{\frac{\alpha-1}{2}}u(x)\cdot
        \textbf{D}^{\frac{\alpha-1}{2}\ast}\psi_1(x)\,dx, \vspace{4pt}\\
        \int_{-1}^1 r(x)\psi_2(x)\,dx
        =\bar{q}~\int_{-1}^1 \textbf{D}^{\frac{\alpha-1}{2}\ast}u(x)\cdot
        \textbf{D}^{\frac{\alpha-1}{2}}\psi_2(x)\,dx,  \vspace{4pt}\\
        \int_{-1}^1 [l(x)+r(x)]D \psi_3(x)\,dx
        -\bar{d}\cdot\int_{-1}^1 u(x)D \psi_3(x)\,dx=\int_{-1}^1 h(x)\psi_3(x)\,dx.
    \end{array}
  \right.
\end{split}
\end{equation}

If $u(x)\in \Phi_2^{\frac{\alpha-1}{2}}(\Omega)$, $\psi_1(x)$ and $\psi_2(x)$ belong to $\hat{I}^{\frac{\alpha-1}{2}}[L_2(\Omega)]$, then $l(x)$ and $r(x)$ belong to $L_{q_2}(-1,1)$, $\psi_3(x)\in W_{0}^{1,p_2}(-1,1)$, and $h(x)\in W^{-1,q_2}(-1,1)$, where $p_2=\frac{2}{2-\alpha}$, $q_2=\frac{2}{\alpha}$. In this case, (\ref{equation_weak_integration3}) is actually the same as Variational Formulation \uppercase\expandafter{\romannumeral2} (\ref{weak_form_2}).

For the convenience of computation and implementation, we partially yield to the requirements on the regularity in (\ref{equation_weak_integration2}). Specifically, we define the third type of mixed variational formulation of (\ref{main_equation}) in the following way:
\begin{definition}[\textrm Variational Formulation \uppercase\expandafter{\romannumeral3}]\label{weak_form_3}
Find $u(x)\in \Phi_2^{\frac{\alpha-1}{2}}(\Omega)$, $l(x)\in L_{2}(-1,1)$, and $r(x)\in L_{2}(-1,1)$, such that
\begin{equation}\label{equation:weak_form_3}
\begin{split}
    \left\{ \begin{array}{ll}
        \left( l(x),\psi_1(x)\right)
        -\bar{p}~\left( \textbf{D}^{\frac{\alpha-1}{2}}u(x),\textbf{D}^{\frac{\alpha-1}{2}\ast}\psi_1(x)\right)=0 & \forall \psi_1(x)\in \hat{I}^{\frac{\alpha-1}{2}}[L_2(-1,1)], \vspace{4pt}\\
        \left( r(x),\psi_2(x)\right)
        -\bar{q}~\left( \textbf{D}^{\frac{\alpha-1}{2}\ast}u(x),\textbf{D}^{\frac{\alpha-1}{2}}\psi_2(x)\right)=0 & \forall \psi_2(x)\in \hat{I}^{\frac{\alpha-1}{2}}[L_2(-1,1)],   \vspace{4pt}\\
        \left( l(x)+r(x),D \psi_3(x)\right)-\bar{d}\cdot\left( u(x),D \psi_3(x)\right)=F_3(\psi_3)
        & \forall \psi_3(x)\in H_0^{1}(-1,1),
    \end{array}
  \right.
\end{split}
\end{equation}
where for a given $h(x)\in H^{-1}(-1,1)$, $F_3:H_0^{1}(-1,1)\rightarrow \mathbb{R}$ is a linear functional defined as
\begin{equation}\label{functional_form_3}
F_3(v):=\langle h,v\rangle,
\end{equation}
and $\langle \cdot,\cdot \rangle$ is the duality pair of $H^{-1}(-1,1)$ and $H_{0}^{1}(-1,1)$.
\end{definition}

We can see that the main difference between the weak formulae (\ref{equation:weak_form_2}) and (\ref{equation:weak_form_3}) is that the linear functional $h(x)$ of later one lies in a bit smaller space than that of former one , and as a sequence, the weak solution $u(x)$ of (\ref{equation:weak_form_3}) lies in a smaller space.


Denote
$$\Psi_{3,N}(\Omega)=\left\{f:f\in P_N(\Omega), \textrm{~and~}f(-1)=f(1)=0\right\}.$$
Then the Galerkin approximation of (\ref{equation:weak_form_3}) is: find $u_{3,N}(x)\in \Phi_{2,N}^{\frac{\alpha-1}{2}}(\Omega)$, $l_N(x)\in P_{N}(\Omega)$, and $r_N(x)\in P_{N}(\Omega)$, such that
\begin{equation}\label{equation:mixed_form_3}
\begin{split}
    \left\{ \begin{array}{ll}
        \left( l_N,\psi_{1,N}\right)
        -\bar{p}~\left( \textbf{D}^{\frac{\alpha-1}{2}}u_{3,N},\textbf{D}^{\frac{\alpha-1}{2}\ast}\psi_{1,N}\right)=0 & \forall \psi_{1,N}\in \hat{I}^{\frac{\alpha-1}{2}}[P_{N}(\Omega)], \vspace{4pt}\\
        \left( r_N,\psi_{2,N}\right)
        -\bar{q}~\left( \textbf{D}^{\frac{\alpha-1}{2}\ast}u_{3,N},\textbf{D}^{\frac{\alpha-1}{2}}\psi_{2,N}\right)=0 & \forall \psi_{2,N}\in \hat{I}^{\frac{\alpha-1}{2}}[P_{N}(\Omega)],   \vspace{4pt}\\
        \left( l_N+r_N,D \psi_{3,N}\right)-\bar{d}\cdot\left( u_{3,N},D \psi_{3,N}\right)=F_3(\psi_{3,N})
        & \forall \psi_{3,N}\in \Psi_{3,N}(\Omega).
    \end{array}
  \right.
\end{split}
\end{equation} 

\subsection{Numerical implementation}\label{subsection:3.2}
In this paper, we mainly focus on designing the numerical schemes for the above variational formulations, and leave the theoretical part in our future work.

We shall make use of the so-called Generalized Jacobi functions that we mentions in Remark \ref{remark1} and have been widely used in other papers of spectral methods for fractional problem, such as \cite{Tian:14,Li:12,Zayernouri:13,Chen:16} and so on.

Recall the following formulas (\cite{Askey:75}, p.20):

\begin{equation}\label{left_jacobi_integral}
\,_{-1}I_{x}^{\alpha}\big((1+x)^{\delta} J_{n}^{\gamma,\delta}(x)\big)
=\frac{\Gamma(n+\delta+1)}{\Gamma(n+\delta+\alpha+1)}
(1+x)^{\delta+\alpha} J_{n}^{\gamma-\alpha,\delta+\alpha}(x),
\end{equation}

\begin{equation}\label{right_jacobi_integral}
\,_{x}I_{1}^{\alpha}\big((1-x)^{\delta} J_{n}^{\delta,\gamma}(x)\big)
=\frac{\Gamma(n+\delta+1)}{\Gamma(n+\delta+\alpha+1)}
(1-x)^{\delta+\alpha} J_{n}^{\delta+\alpha,\gamma-\alpha}(x),
\end{equation}
where $\alpha>0$, $\delta>-1$, $\gamma\in\mathbb{R}$.

Using the properties
$$\textbf{D}^{\alpha}\,_{-1}I_{x}^{\alpha}=I,$$
and
$$\textbf{D}^{\alpha\ast}\,_{x}I_{1}^{\alpha}=I,$$
we can get from formulae (\ref{left_jacobi_integral}) and
(\ref{right_jacobi_integral}) respectively that
\begin{equation}\label{left_jacobi_derivative}
\textbf{D}^{\alpha}\big((1+x)^{\delta+\alpha} J_{n}^{\gamma-\alpha,\delta+\alpha}(x)\big)
=\frac{\Gamma(n+\delta+\alpha+1)}{\Gamma(n+\delta+1)}(1+x)^{\delta} J_{n}^{\gamma,\delta}(x),
\end{equation}
\begin{equation}\label{right_jacobi_derivative}
\textbf{D}^{\alpha\ast}\big((1-x)^{\delta+\alpha} J_{n}^{\delta+\alpha,\gamma-\alpha}(x)\big)
=\frac{\Gamma(n+\delta+\alpha+1)}{\Gamma(n+\delta+1)}(1-x)^{\delta} J_{n}^{\delta,\gamma}(x).
\end{equation}

\subsubsection{Galerkin spectral scheme of Variational Formulation-\uppercase\expandafter{\romannumeral1}}
\label{subsection:3.2.1}

For the discrete variational formulation (\ref{equation:Galerkin_form_1}), we construct two kinds of trial functions as
\begin{equation}\label{left_trial_1}
\phi^L_{1,n}(x):={}_{-1}I_{x}^{\frac{\alpha}{2}} L_n(x)
=\frac{\Gamma(n+1)}{\Gamma(n+1+\frac{\alpha}{2})}
(1+x)^{\frac{\alpha}{2}}J_n^{-\frac{\alpha}{2},\frac{\alpha}{2}}(x),
~0\leq n\leq N-1,
\end{equation}
and
\begin{equation}\label{right_trail_1}
\phi^R_{1,n}(x):={}_{x}I_{1}^{\frac{\alpha}{2}} L_n(x)
=\frac{\Gamma(n+1)}{\Gamma(n+1+\frac{\alpha}{2})}
(1-x)^{\frac{\alpha}{2}}J_n^{\frac{\alpha}{2},-\frac{\alpha}{2}}(x),
~0\leq n\leq N-1,
\end{equation}
where $L_n(x)=J_{n}^{0,0}(x),~n\geq0$, are Legendre polynomials, which are orthogonal in the $L_2$
sense \cite{Canuto:06,Hesthaven:07}:
\begin{equation*}
\int_{-1}^{1}L_n(x)L_m(x)=\gamma_n\delta_{mn},~~\gamma_n=\frac{2}{2n+1};
\end{equation*}
and take test functions as
\begin{equation}\label{test_1}
  v_{1,k}(x):=(1+x){}_xI_1^{\frac{\alpha}{2}}L_k(x), \qquad 0\leq k\leq N-1.
\end{equation}

Denote
$$u_{1,N}(x):=\sum_{n=0}^{N-1} u_{1,n}^L \phi^L_{1,n}(x)$$
be the approximation of the exact solution $u$, and let
\begin{equation}\label{numerical_solution_1}
u_{1,N}(x_i)
=\sum_{n=0}^{N-1} u_{1,n}^L \phi^L_{1,n}(x_i)
= \sum_{m=0}^{N-1} u_{1,m}^R \phi^R_{1,m}(x_i),~i=1,\cdots,N,
\end{equation}
for some given nodes $\{x_i\}_{i=1}^{N}$.

Denote
$$ \textbf{u}_1^L=[u^L_{1,0},u^L_{1,1},\cdots,u^L_{1,N-1}]^T,\quad
\textbf{u}_1^R=[u^R_{1,0},u^R_{1,1},\cdots,u^R_{1,N-1}]^T, $$
and
$A_1^{L}$ and $A_1^{R}$ as two $N\times N$ matrices with
$$(A_1^L)_{i,j}=\phi^L_{1,j-1}(x_i), \quad (A_1^R)_{i,j}=\phi^R_{1,j-1}(x_i).$$
Then \eqref{numerical_solution_1} can be rewritten as
\begin{equation}\label{u_right_to_left_1}
  A_1^L\, \textbf{u}_1^L= A_1^R\, \textbf{u}_1^R.
\end{equation}

Use the properties of Legendre polynomials \cite{Hesthaven:07}
\begin{equation}\label{legendre_property_1}
(2k+1)L_k(x)=\frac{d}{dx}(L_{k+1}(x)-L_{k-1}(x)),
\end{equation}
\begin{equation}\label{legendre_property_2}
L_k(\pm 1)=(\pm1)^k,
\end{equation}
and Leibniz rule for fractional derivative \cite{Podlubny:99}, we can obtain from (\ref{test_1}) that
$$\textbf{D}^{\frac{\alpha}{2}\ast}v_{1,k}(x)
=\left\{\begin{array}{lll}
(1+x)L_k(x)+\frac{\alpha}{2(2k+1)}(L_{k+1}(x)-L_{k-1}(x)),& k\geq 1,\vspace{5pt}\\
(1+x)L_k(x)+\frac{\alpha}{2(2k+1)}(L_{k+1}(x)-1),& k=0.
\end{array}\right.
$$

For computing the left fractional derivative of $v_{1,k}(x)$, we denote
\begin{equation}\label{test_left_1}
  v_{1,k}(x)=\sum_{n=0}^\infty (v_{1,k})_{n}~{}_{-1}I_x^{\frac{\alpha}{2}} L_n(x).
\end{equation}
Taking inner product with $\textbf{D}^{\frac{\alpha}{2}\ast} L_m(x)$ in \eqref{test_left_1}, using the orthogonality of Legendre polynomials and formulae \eqref{right_jacobi_integral}, \eqref{right_jacobi_derivative}, one obtains
$$
(v_{1,k})_{m}
=\left(m+\frac{1}{2}\right)\frac{\Gamma(k+1)\Gamma(m+1)}{\Gamma(k+1+\alpha/2)\Gamma(m+1-\alpha/2)}
\Big((1+x)J_k^{\frac{\alpha}{2},-\frac{\alpha}{2}}(x),J_m^{-\frac{\alpha}{2},\frac{\alpha}{2}}(x)\Big).$$

Therefore, the matrix formulation of \eqref{equation:Galerkin_form_1} is
\begin{equation}\label{matrix_form_left_right_1}
  -\bar{p}\cdot M_1^L \,\textbf{u}_1^L-\bar{q}\cdot M_1^R \,\textbf{u}_1^R+\bar{d}\cdot \,M_1^C \,\textbf{u}_1^L=\textbf{f}_1,
\end{equation}
where
\begin{equation*}
  \begin{array}{lll}
   &&(M_1^L)_{k+1,n+1}\vspace{5pt}\\
   &=&\left(L_n(x),(1+x)L_k(x)+\frac{\alpha}{2(2k+1)}\left(L_{k+1}(x)-L_{k-1}(x)\right)\right)\vspace{4pt}\\
   &=&\left(L_n(x),xL_k(x)\right)+\gamma_k\delta_{n,k}
   +\frac{\alpha}{2(2k+1)}\left(\gamma_{k+1}\delta_{n,k+1}-\gamma_{k-1}\delta_{n,k-1}\right),
   \quad k\geq 1;\vspace{10pt}\\
   &&(M_1^L)_{1,n+1}\vspace{5pt}\\
   &=&\left(L_n(x),(1+x)L_0(x)+\frac{\alpha}{2}\left(L_{1}(x)-1\right)\right)\vspace{4pt}\\
   &=&(2-\alpha)\delta_{n,0}+\frac{2+\alpha}{3}\delta_{n,1};\vspace{10pt}\\
   &&(M_1^R)_{k+1,n+1}\vspace{5pt}\\
   &=&
   (L_n(x),\sum_{m=0}^\infty (v_{1,k})_{m}L_m(x))=(v_{1,k})_{n}\cdot\gamma_n,\vspace{10pt}\\
   \end{array}
\end{equation*}
\begin{equation*}
  \begin{array}{lll}
 &&(M_1^C)_{k+1,n+1}\vspace{5pt}\\
   &=&\Big(\textbf{D}^{\frac{1}{2}}\phi^L_{1,n}(x),\textbf{D}^{\frac{1}{2}\ast}v_{1,k}(x)\Big)\vspace{5pt}\\
   &=& \frac{\Gamma(n+1)\Gamma(k+1)}{\Gamma(n+\beta+1)\Gamma(k+\beta+1)}\cdot\vspace{5pt}\\
       &&\Big((1-x^2)^\beta\cdot J_n^{-\beta,\beta}(x),(1+x) J_k^{\beta,-\beta}(x)-\frac{1-x}{2(k+\beta+1)} J_k^{1+\beta,-1-\beta}(x)\Big),
   \end{array}
\end{equation*}
and
$$(\textbf{f}_{1}){k}=\frac{\Gamma(k+1)}{\Gamma(k+1+\frac{\alpha}{2})}
\int_{-1}^1 (1-x)^{\frac{\alpha}{2}}(1+x)h(x)J_k^{\frac{\alpha}{2},-\frac{\alpha}{2}}(x)\,dx,$$
with $\beta=\frac{\alpha}{2}$, and $M^C_1$ is calculated by using Leibniz rule for fractional derivative \cite{Podlubny:99}.

It should be noted that all of these integrals in above formulations can be computed exactly by Gauss quadrature or weighted Gauss quadrature.

Combined (\ref{u_right_to_left_1}) with (\ref{matrix_form_left_right_1}), we can get the final Galerkin spectral scheme of \eqref{equation:Galerkin_form_1}:
\begin{equation}\label{scheme_1}
  \bigg( -\bar{p}\cdot M_1^L- \bar{q}\cdot M_{1}^R\,(A_1^R)^{-1}\,A_1^L+\bar{d}\cdot M_1^C \bigg) \textbf{u}_1^L=\textbf{f}_1.
\end{equation}

\subsubsection{Petrov-Galerkin spectral scheme of Variational Formulation-\uppercase\expandafter{\romannumeral2}}
\label{subsection:3.2.2}

For the discrete variational formulation (\ref{equation:Petrov_Galerkin_form_2}), we construct the corresponding two kinds of trial functions as
\begin{equation}\label{left_trial_2}
\phi^L_{2,n}(x):={}_{-1}I_{x}^{\frac{\alpha-1}{2}} L_n(x),
\qquad 0\leq n\leq N-1,
\end{equation}
and
\begin{equation}\label{right_trail_2}
\phi^R_{2,n}(x):={}_{x}I_{1}^{\frac{\alpha-1}{2}} L_n(x),
\qquad 0\leq n\leq N-1,
\end{equation}
and take the corresponding test functions as
\begin{equation}\label{test_2}
  v_{2,k}(x):=(1+x){}_xI_1^{\frac{\alpha+1}{2}}L_k(x), \qquad 0\leq k\leq N-1.
\end{equation}

Denote
$$u_{2,N}(x):=\sum_{n=0}^{N-1} u_{2,n}^L \phi^L_{2,n}(x)$$
as the approximation of the exact solution $u$, and let
\begin{equation}\label{numerical_solution_2}
u_{2,N}(x_i)
=\sum_{n=0}^{N-1} u_{2,n}^L \phi^L_{2,n}(x_i)
= \sum_{m=0}^{N-1} u_{2,m}^R \phi^R_{2,m}(x_i),~i=1,\cdots,N,
\end{equation}
for some given nodes $\{x_i\}_{i=1}^{N}$, and denote
$$ \textbf{u}_2^L=[u^L_{2,0},u^L_{2,1},\cdots,u^L_{2,N-1}]^T,\quad
\textbf{u}_2^R=[u^R_{2,0},u^R_{2,1},\cdots,u^R_{2,N-1}]^T. $$
Similarly, there is
\begin{equation}\label{u_right_to_left_2}
  A_2^L\, \textbf{u}_2^L= A_2^R\, \textbf{u}_2^R,
\end{equation}
where
$A_2^{L}$ and $A_2^{R}$ are two $N\times N$ matrices with
$$(A_2^L)_{i,j}=\phi^L_{2,j-1}(x_i), \quad (A_2^R)_{i,j}=\phi^R_{2,j-1}(x_i).$$

Denote
\begin{equation}\label{test_left_2}
  v_{2,k}(x)=\sum_{n=0}^\infty (v_{2,k})_{n}~{}_{-1}I_x^{\frac{\alpha+1}{2}} L_n(x).
\end{equation}
Again, by using Leibniz rule for fractional derivative \cite{Podlubny:99}, we can get
$$
(v_{2,k})_{m}
=\left(m+\frac{1}{2}\right)\frac{\Gamma(k+1)\Gamma(m+1)}
{\Gamma(k+\frac{\alpha+3}{2})\Gamma(m-\frac{\alpha-1}{2})}
\Big((1+x)J_k^{\frac{\alpha+1}{2},-\frac{\alpha+1}{2}}(x),J_m^{-\frac{\alpha+1}{2},\frac{\alpha+1}{2}}(x)\Big),$$
and
$$\textbf{D}^{\frac{\alpha+1}{2}\ast}v_{2,k}(x)
=\left\{\begin{array}{ll}
(1+x)L_k(x)+\frac{\alpha+1}{2(2k+1)}(L_{k+1}(x)-L_{k-1}(x)),&k\geq 1,\vspace{5pt}\\
(1+x)L_k(x)+\frac{\alpha+1}{2(2k+1)}(L_{k+1}(x)-1),&k=0.
\end{array}\right.$$

Thus, the Petrov-Galerkin spectral scheme of \eqref{equation:Petrov_Galerkin_form_2} is
\begin{equation}\label{scheme_2}
  \left(-\bar{p}\cdot M_2^L -\bar{q}\cdot M_2^R \,(A_2^R)^{-1}\,(A_2^L)+\bar{d}\cdot \,M_2^C\right)\,\textbf{u}_2^L=\textbf{f}_2,
\end{equation}
where
\begin{equation*}
  \begin{array}{lll}
   &&(M_2^L)_{k+1,n+1}\vspace{5pt}\\
   &=&\left\{\begin{array}{lll}
   \left(L_n(x), (1+x)L_k(x)+\frac{\alpha+1}{2(2k+1)}(L_{k+1}(x)-L_{k-1}(x))\right),&k\geq 1,\vspace{5pt}\\
   (1-\alpha)\delta_{n,0}+\frac{\alpha+3}{3}\delta_{n,1},&k=0;
   \end{array}\right.\vspace{10pt}\\
   &&(M_2^R)_{k+1,n+1}\vspace{5pt}\\
   &=&(L_n(x),\sum_{m=0}^\infty (v_{2,k})_{m}L_m(x))
   =(v_{2,k})_{n}\cdot\gamma_n,\vspace{5pt}\\
   \end{array}
\end{equation*}
\begin{equation*}
   \begin{array}{lll}
   &&(M_2^C)_{k+1,n+1}\vspace{5pt}\\
   &=&\Big(\phi^L_{2,n}(x),-D\,v_{2,k}(x)\Big)\vspace{5pt}\\
   &=& \frac{\Gamma(n+1)\Gamma(k+1)}{\Gamma(n+\beta+1)\Gamma(k+\beta+1)}\cdot  \vspace{5pt}\\
    && \Big((1-x^2)^\beta J_n^{-\beta,\beta}(x),
       (1+x) J_k^{\beta,-\beta}(x)-\frac{1-x}{(k+\beta+1)} J_k^{1+\beta,-1-\beta}(x)\Big),
   \end{array}
\end{equation*}
with $\beta=\frac{\alpha-1}{2}$, and
$$(\textbf{f}_{2})_{k}=\frac{\Gamma(k+1)}{\Gamma(k+\frac{\alpha+3}{2})}
\int_{-1}^1 (1-x)^{\frac{\alpha+1}{2}}(1+x)h(x)J_k^{\frac{\alpha+1}{2},-\frac{\alpha+1}{2}}(x)\,dx.$$

During the computation, it is found that the condition numbers of the stiffness matrix in the above two schemes (\ref{scheme_1}) and (\ref{scheme_2}) are increasing as $O(N^{2\alpha})$. When $\alpha$ is close to $2$, the condition numbers increase fast, making the numerical solution sensitive to a small disturbance. The usual method to deal with ill-conditioned system is precondition. However, the stiffness matrices here are full, which makes it difficult to find an appropriate preconditioning matrix for them. While, at the cost of losing a bit of regularity for the solution, the mixed Galerkin spectral system introduced below, instead, shows to be well-conditioned.

\subsubsection{Mixed Galerkin spectral scheme of Variational Formulation-\uppercase\expandafter{\romannumeral3}}
\label{subsection:3.2.3}

For the discrete variational formulation (\ref{equation:mixed_form_3}), we obtain the matrix form by four steps.

Step 1: Similar to the previous two schemes, we construct two kinds of trial functions with left and right fractional integrals with $\beta=\frac{\alpha-1}{2}$:
\begin{equation*}
\phi^L_{3,n}(x):={}_{-1}I_{x}^{\beta} L_n(x) \qquad 0\leq n\leq N-1,
\end{equation*}
\begin{equation*}
\phi^R_{3,n}(x):={}_{x}I_{1}^{\beta} L_n(x) \qquad 0\leq n\leq N-1.
\end{equation*}
Denote
\begin{equation}\label{basis-u3}
  u_{3,N}(x):=\sum_{n=0}^{N-1}u_{3,n}^L\phi_{3,n}^L(x)
\end{equation}
as the approximation of the exact solution $u$, and let
$$u_{3,N}(x_i)=\sum_{n=0}^{N-1} u_{3,n}^L\phi_{3,n}^L(x_i)=\sum_{m=0}^{N-1} u_{3,m}^R\phi^R_{3,m}(x_i), ~~i=1,\cdots,N,$$
for some given nodes $\{x_i\}_{i=1}^N$, and denote
\begin{equation}
  \uu_3^L=[u_{3,0}^L,u_{3,1}^L,\cdots,u_{3,N-1}^L]^T, \quad  \uu_3^R=[u_{3,0}^R,u_{3,1}^R,\cdots,u_{3,N-1}^R]^T.
\end{equation}
Similarly, there is
\begin{equation}\label{mix_mat1}
  A_3^L\uu_3^L=A_3^R\uu_3^R,
\end{equation}
where
$A_3^L=A_2^L$ and $A_3^R=A_2^R$.

Step 2: We deal with the first equation of \eqref{equation:mixed_form_3}. Take the trial functions of $l_N(x)$ to be $L_n(x)+L_{n+1}(x)$, and denote
\begin{equation*}
  l_N(x)=\sum_{n=0}^{N-1}l_n (L_n(x)+L_{n+1}(x)).
\end{equation*}
Take the test function $\psi_{1,N}(x)$ to be ${}_xI_1^{\beta}L_k(x),0\leq k\leq N-1$. Substituting them into \eqref{equation:mixed_form_3}, we have
$$\sum_{n=0}^{N-1}l_n(L_n(x)+L_{n+1}(x),{}_xI_1^{\beta}L_k(x))
=\bar{p}\,\sum_{n=0}^{N-1}u_{3,n}^L(L_n(x),L_k(x)),$$
i.e.,
\begin{equation}\label{mix_mat2}
  L^{(\beta)}\cdot \textbf{{l}}= \bar{p}\,B\cdot \uu^L_3,
\end{equation}
where
$$L^{(\beta)}_{k+1,n+1}=(L_n(x)+L_{n+1}(x),{}_xI_1^{\beta}L_k(x)),\quad \textbf{{l}}=[l_0,l_1,\cdots,l_{N-1}]^T,\quad B=\textrm{diag}(\gamma_k).$$

Step 3: We deal with the second equation of \eqref{equation:mixed_form_3}. Take the trial functions of $r_N(x)$ to be $(L_n(x)-L_{n+1}(x))$, and denote $$r_N(x)=\sum_{n=0}^{N-1}r_n (L_n(x)-L_{n+1}(x)).$$
Take the test function $\psi_{2,N}(x)$ to be ${}_{-1}I_x^{\beta}L_k(x),0\leq k\leq N-1$.
Substituting them into \eqref{equation:mixed_form_3}, we have
$$\sum_{n=0}^{N-1}r_n(L_n(x)-L_{n+1}(x),{}_{-1}I_x^{\beta}L_k(x))
=\bar{q}\,\sum_{n=0}^{N-1}u_{3,n}^R(L_n(x),L_k(x)),$$
i.e.,
\begin{equation}\label{mix_mat3}
  R^{(\beta)}\cdot \textbf{{r}}= \bar{q}\,B\cdot \uu_3^R,
\end{equation}
where
$$R^{(\beta)}_{k+1,n+1}=(L_n(x)-L_{n+1}(x),{}_{-1}I_x^{\beta}L_k(x)),\quad \textbf{{r}}=[r_0,r_1,\cdots,r_{N-1}]^T.$$

Step 4: Finally we deal with the third equation of \eqref{equation:mixed_form_3}. Considering the condition $\psi_{3,N}(x)\in H_0^1(\Omega)$, we take the test function to be $L_{k-1}(x)-L_{k+1}(x),1\leq k\leq N$. Substituting them into \eqref{equation:mixed_form_3}, we have
\begin{equation}\label{mix_mat4}
  -C^L\cdot \textbf{{l}} - C^R\cdot \textbf{{r}}  +d\, M_3^C \uu^L=\f_3,
\end{equation}
where
\begin{equation*}
  \begin{split}
  &C^L_{k,n+1}=\Big((L_n(x)+L_{n+1}(x)),(2k+1)L_k(x)\Big)=2(\delta_{k,n}+\delta_{k,n+1}),\\
  &C^R_{k,n+1}=\Big((L_n(x)-L_{n+1}(x)),(2k+1)L_k(x)\Big)=2(\delta_{k,n}-\delta_{k,n+1}),\\
  &(M_3^C)_{k,n+1}=(\phi_n^L(x),(2k+1)L_k(x)),\\
  &(\f_3)_k=(h(x),L_{k-1}(x)-L_{k+1}(x)).
  \end{split}
\end{equation*}

Combining \eqref{mix_mat1}, \eqref{mix_mat2}, \eqref{mix_mat3}, and \eqref{mix_mat4}, we obtain
the mixed Galerkin spectral scheme of (\ref{equation:weak_form_3}) as
\begin{equation}
  (-\bar{p}M_3^L - \bar{q}M_3^R(A^R)^{-1} A^L +\bar{d} M_3^C) \uu^L = \f_3
\end{equation}
with
\begin{equation*}
  M_3^L:=C^L (L^{(\beta)})^{-1} B,  \quad M_3^R:=C^R (R^{(\beta)})^{-1} B.
\end{equation*}

\begin{remark}
  Although we use four steps in the mixed Galerkin spectral scheme which seems a bit more complicated, many matrices in this scheme are sparse and the elements are more convenient to be calculated than those in the previous two schemes.
\end{remark}

\begin{remark}
  Here we make an rough explanation about the reason why we choose the basis functions of $l_N(x)$ as $L_n(x)+L_{n+1}(x)$ which vanishing at $x=-1$. Based on the first equation of  \eqref{equation:mixed_form_3}, we find the left-hand side can be reformed as
  \begin{equation}
    \begin{split}
      (l_N,\psi_{1,N})
      =(l_N,{}_xI_b^{\frac{\alpha-1}{2}}\textbf{D}^{\frac{\alpha-1}{2}\ast} \psi_{1,N} )
      =({}_aI_x^{\frac{\alpha-1}{2}}l_N,\textbf{D}^{\frac{\alpha-1}{2}\ast} \psi_{1,N} ),
    \end{split}
  \end{equation}
 which means
  \begin{equation}
    {}_aI_x^{\frac{\alpha-1}{2}}l_N = \textbf{D}^{\frac{\alpha-1}{2}} u_{3,N},
  \end{equation}
  in the $L_2$ sense. If $l_N(x)$ does not tend to zero as $x\rightarrow-1$, then
  \begin{equation*}
    \textbf{D}^{\frac{\alpha-1}{2}} u_{3,N} \rightarrow O(1+x)^{\frac{\alpha-1}{2}}, \quad x\rightarrow-1;
  \end{equation*}
  On the contrary, the basis functions of $u_{3,N}$ in \eqref{basis-u3} implies
  \begin{equation*}
    \textbf{D}^{\frac{\alpha-1}{2}} u_{3,N} \rightarrow O(1+x)^{\frac{3-\alpha}{2}}, \quad x\rightarrow-1;
  \end{equation*}
which is contradict to each other. Therefore, we restrict $l_N(x)$ to be zero at $x=-1$. Similarly, the basis functions of $r_N(x)$ are chosen as $L_n(x)-L_{n+1}(x)$ that vanishing at $x=1$.
\end{remark}

\section{Numerical tests}\label{section:4}

In what follows, we provide some numerical results to verify the validity of our proposed three kinds of numerical schemes---Galerkin spectral scheme, Petrov-Galerkin spectral scheme and, most impotently, mixed Galerkin spectral scheme. For convenience, we denote them as Scheme 1, Scheme 2 and Scheme 3, respectively, in this section.

\begin{figure}[ht]
\begin{minipage}{0.45\linewidth}
  \centerline{\includegraphics[scale=0.4]{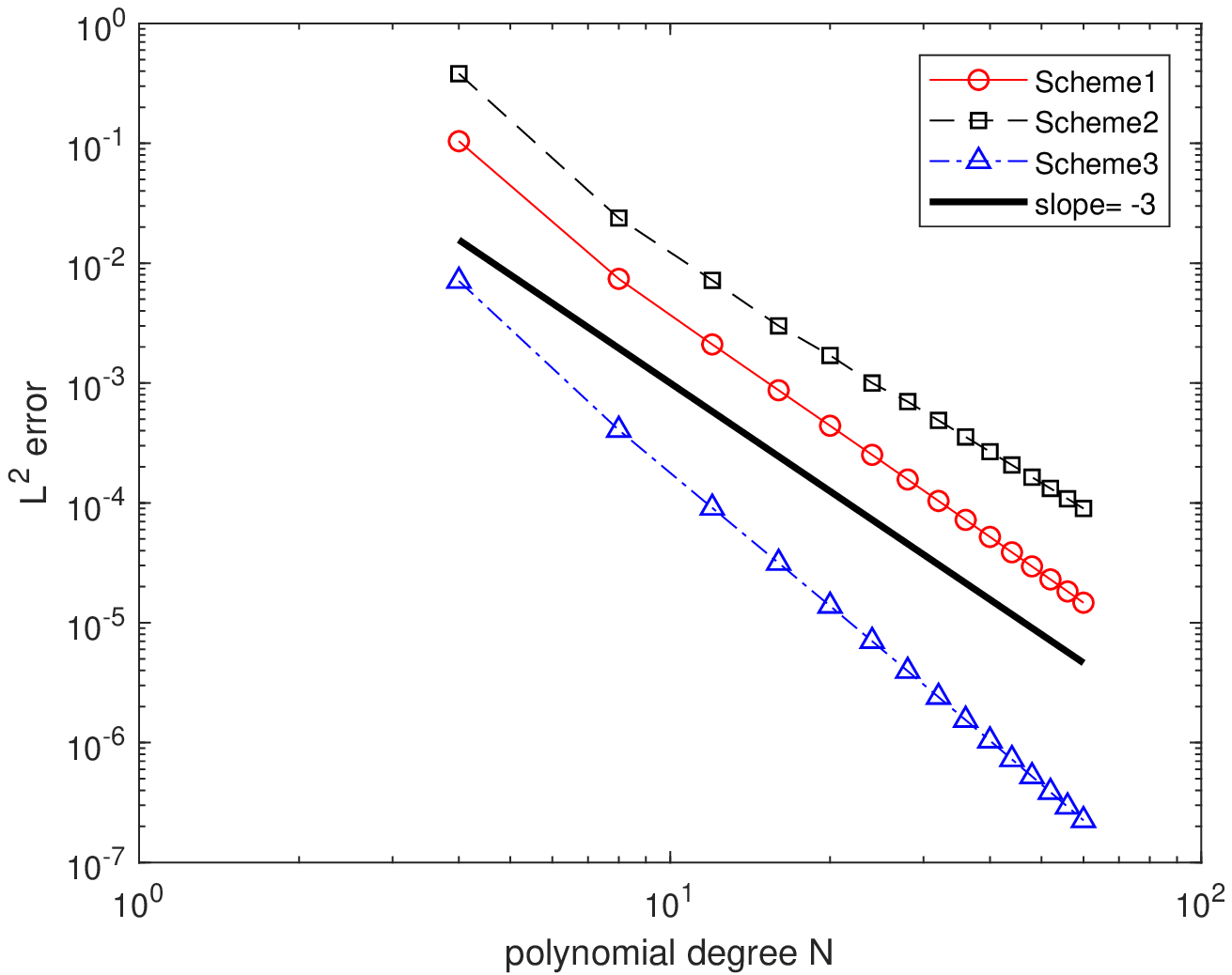}}
  \centerline{(a)$\alpha=1.3$}
\end{minipage}
\hfill
\begin{minipage}{0.45\linewidth}
  \centerline{\includegraphics[scale=0.4]{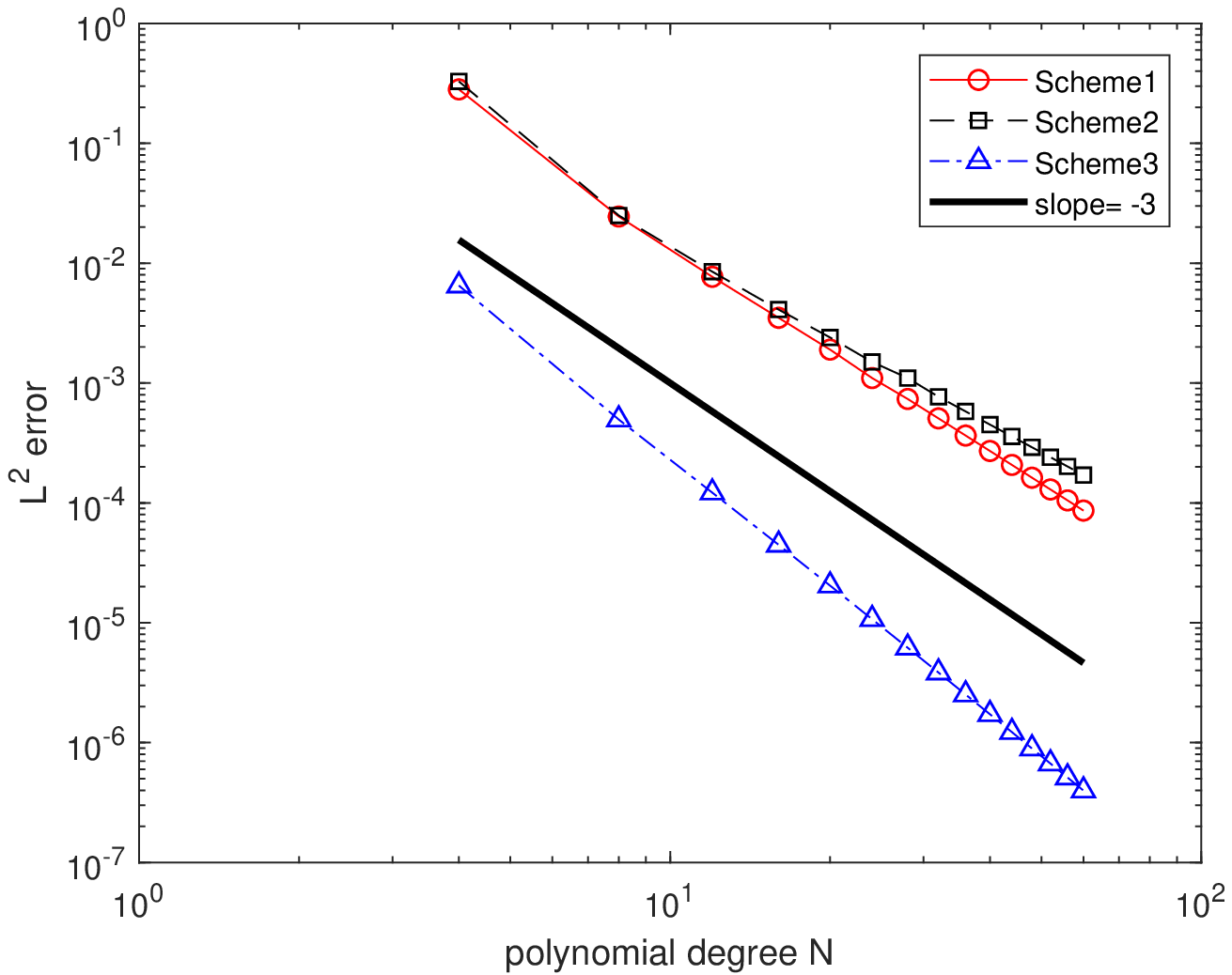}}
  \centerline{(b) $\alpha=1.6$}
\end{minipage}
\caption{The numerical $L_2$ errors of Example \ref{exam1} vs the polynomial degree N for three different schemes with $\alpha=1.3$ in (a) and $\alpha=1.6$ in (b).}
\label{fig1}
\end{figure}

\begin{example}\label{exam1}
We firstly consider a one-sided problem to verify the effectiveness of our numerical schemes. More precisely, let
$\bar{p}=1$ and $\bar{q}=\bar{d}=0$ in \eqref{main_equation}, i.e.,
  \begin{equation}
    {}_{-1}D_x^{\alpha}u(x)=h(x).
  \end{equation}
  We choose the exact solution to be
  \begin{equation}
    u(x)=\left\{
    \begin{array}{ll}
      -\frac{(x+1)^{3+\frac{\alpha}{2}}}{\Gamma(4+\frac{\alpha}{2})}   &~ [-1,0], \\[5pt]
      \frac{2x^{3+\frac{\alpha}{2}}-(x+1)^{3+\frac{\alpha}{2}}}{\Gamma(4+\frac{\alpha}{2})}   &~(0,1],
    \end{array}\right.
  \end{equation}
  so that the source term is
  \begin{equation}
    h(x)=\left\{
    \begin{array}{ll}
      -\frac{(x+1)^{3-\frac{\alpha}{2}}}{\Gamma(4-\frac{\alpha}{2})}   &~ [-1,0], \\[5pt]
      \frac{2x^{3-\frac{\alpha}{2}}-(x+1)^{3-\frac{\alpha}{2}}}{\Gamma(4-\frac{\alpha}{2})}   &~(0,1].
    \end{array}\right.
  \end{equation}
\end{example}
After applying $\frac{\alpha}{2}$-order fractional derivative, the solution has a finite regularity at the point $x=0$, i.e.,
\begin{equation}
    v(x)={}_{-1}D_x^{\frac{\alpha}{2}}u(x)=\left\{
    \begin{array}{ll}
      -\frac{(x+1)^3}{\Gamma(4)}   &~ [-1,0], \\[5pt]
      \frac{2x^{3}-(x+1)^{3}}{\Gamma(4)}   &~(0,1].
    \end{array}\right.
\end{equation}
In fact, $v(x)\in H^{3+\frac{1}{2}-\epsilon}[-1,1]$, for $\epsilon>0$, since
    \begin{equation}
    D^3v(x)=\left\{
    \begin{array}{ll}
      -1   &~ [-1,0], \\[5pt]
      1   &~(0,1].
    \end{array}\right.
\end{equation}
We plot the $L_2$ error vs the polynomial degree $N$ for $\alpha=1.3$ and $\alpha=1.6$ in Figure \ref{fig1}, and find that the convergence order is around $N^{-3.5}$. Although the convergence orders of the three schemes look similar, the magnitude of the error of Scheme 3 is smaller than the other ones.

\begin{figure}[ht]
\begin{minipage}{0.31\linewidth}
  \centerline{\includegraphics[scale=0.3]{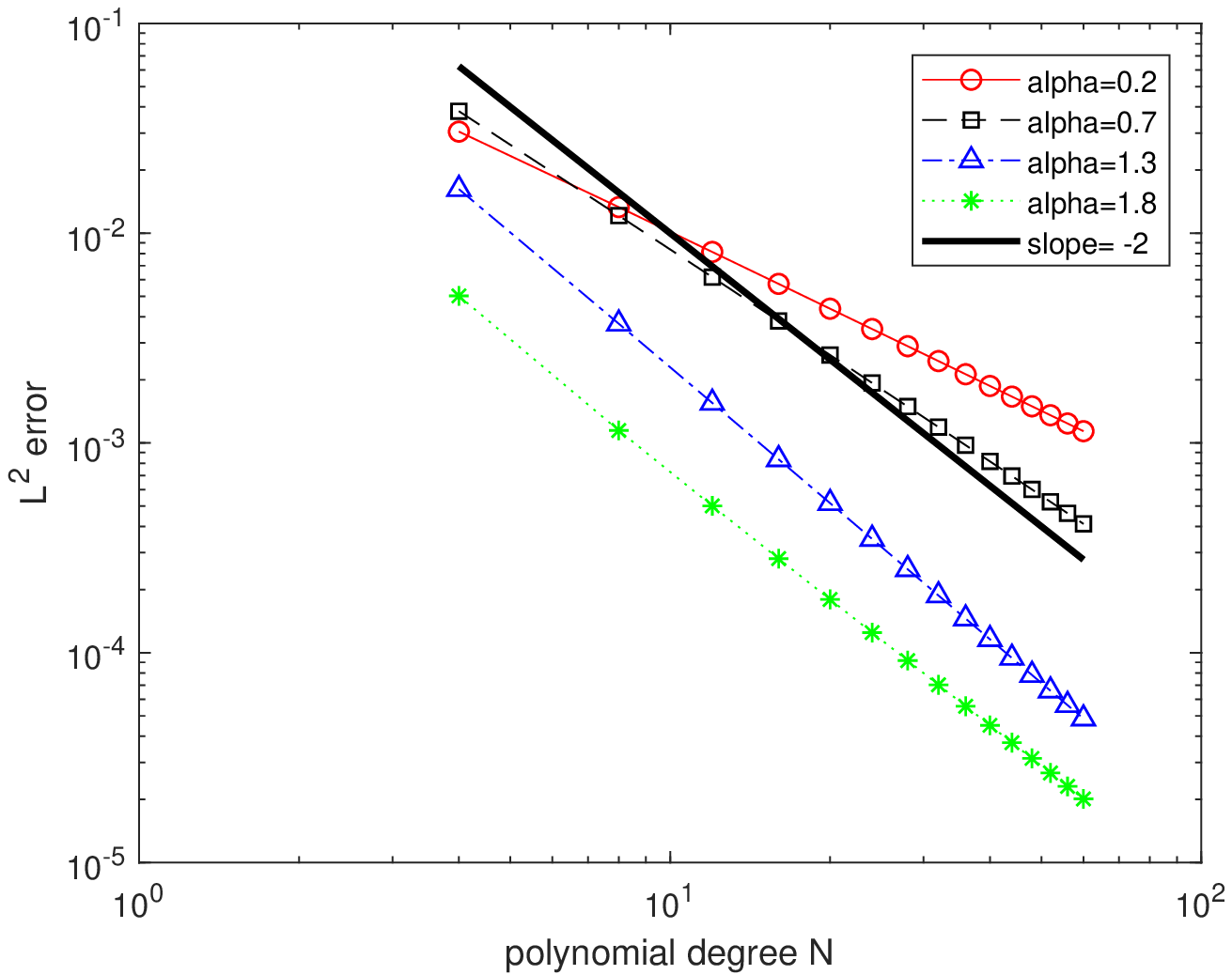}}
  \centerline{(a) Scheme 1}
\end{minipage}
\hfill
\begin{minipage}{0.31\linewidth}
  \centerline{\includegraphics[scale=0.3]{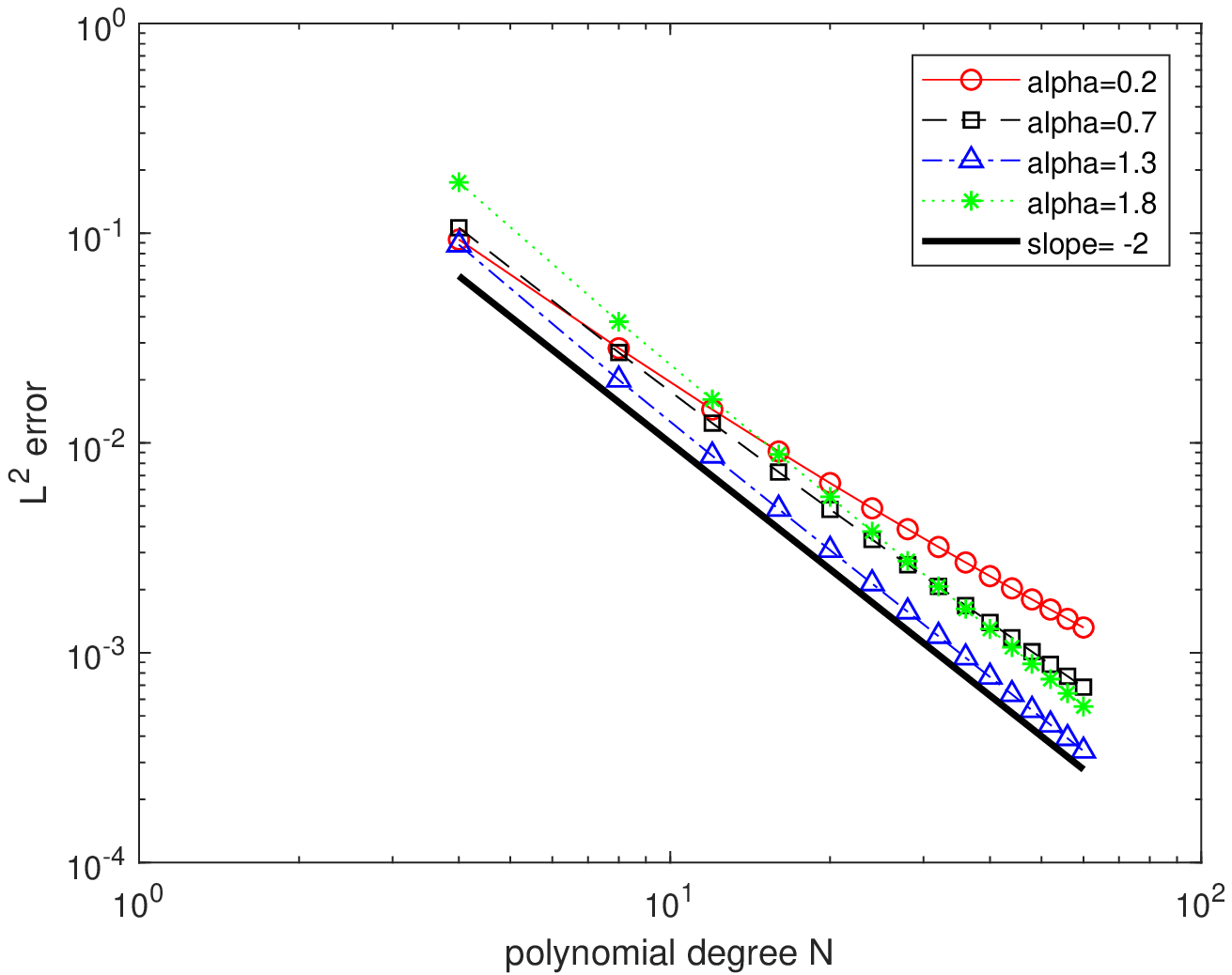}}
  \centerline{(b) Scheme 2}
\end{minipage}
\hfill
\begin{minipage}{0.31\linewidth}
  \centerline{\includegraphics[scale=0.3]{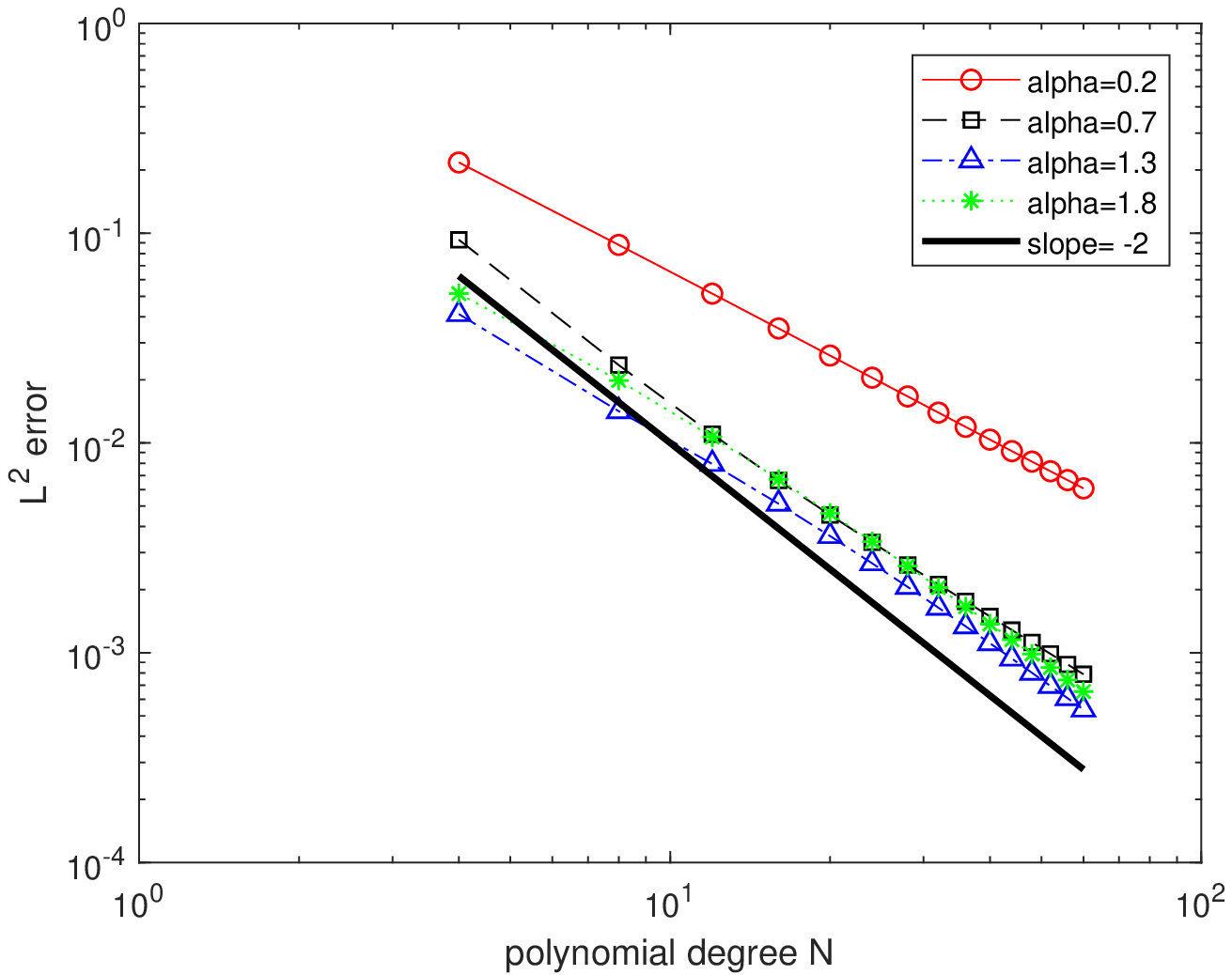}}
  \centerline{(c) Scheme 3}
\end{minipage}
\caption{The numerical $L_2$ errors of Example \ref{exam2} vs the polynomial degree N with three different schemes in
(a), (b) and (c), respectively.}
\label{fig2}
\end{figure}

\begin{example}\label{exam2}
In this example, we apply the three schemes to the fractional Laplacian equation in one dimension case, i.e.,
  \begin{equation}
  \left\{
    \begin{array}{rll}
    (-\Delta)^{\alpha/2}u(x)&=~h(x), &  x\in(-1,1) \\[5pt]
    u(x)&=~0, & x\in \mathbb{R}\backslash(-1,1).
    \end{array}\right.
  \end{equation}
For the source term $h(x)=1$ in $(-1,1)$, the exact solution is \cite{Getoor:61}
 \begin{equation}\label{solution-Laplace}
   u(x)=\frac{2^{-\alpha}\Gamma(\frac{1}{2})}{\Gamma(\frac{1+\alpha}{2})\Gamma(1+\frac{\alpha}{2})}
 (1-x^2)^{\alpha/2} \quad \textrm{in~}(-1,1).
 \end{equation}
\end{example}
The numerical tests for this fractional Laplacian equation are presented in Figure \ref{fig2}. We find that the errors all decay algebraically (about O($N^{-2}$)), which implies our proposed three schemes are all valid not only for $\alpha\in(1,2)$, but also for the whole range of $\alpha\in(0,2)$. In addition, we can see that even for the solution $u(x)$ with low regularity such as in \eqref{solution-Laplace}, all of the three schemes show a good convergence result.

\begin{figure}[ht]
\begin{minipage}{0.45\linewidth}
  \centerline{\includegraphics[scale=0.4]{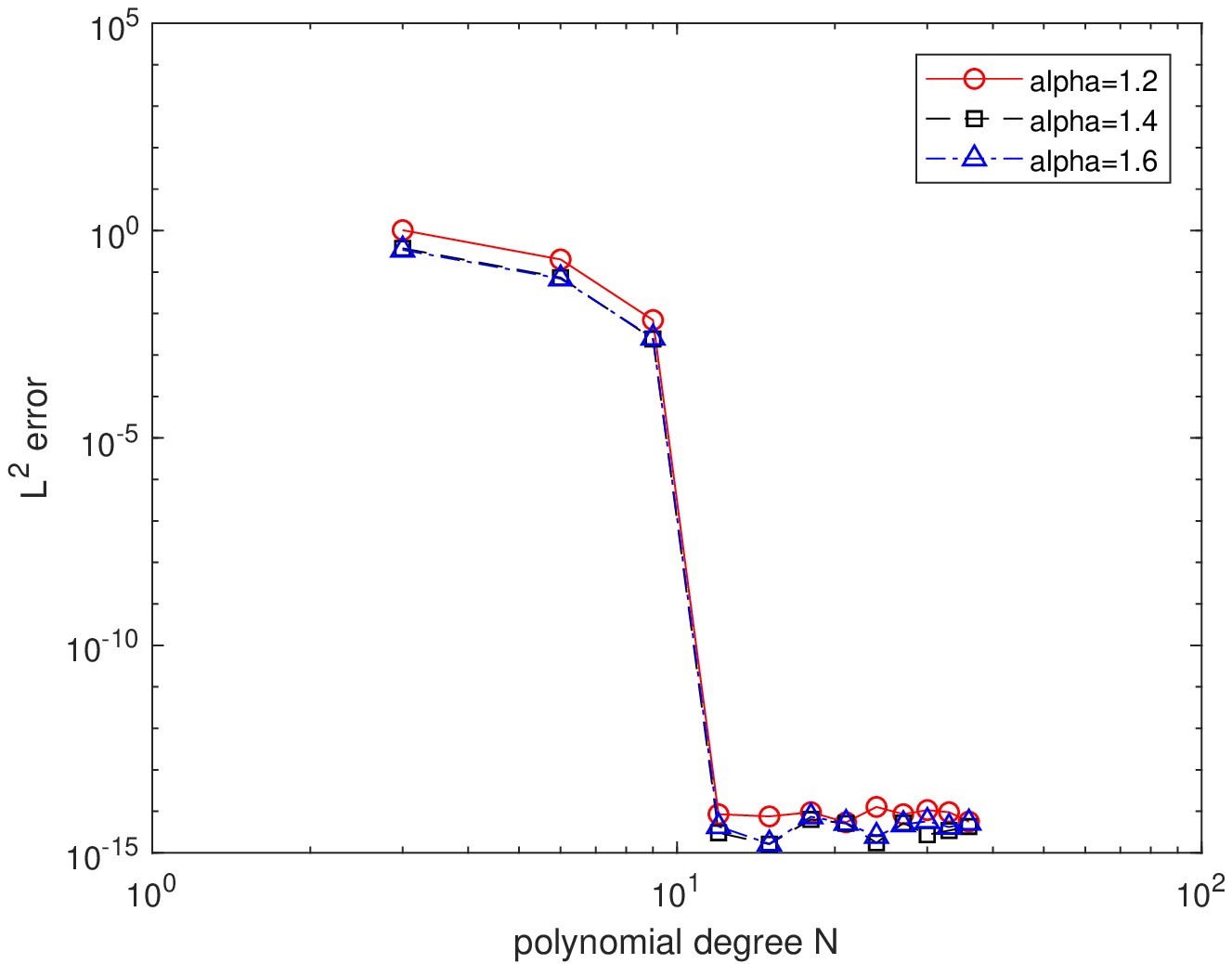}}
  \centerline{(a) Scheme 1}
\end{minipage}
\hfill
\begin{minipage}{0.45\linewidth}
  \centerline{\includegraphics[scale=0.4]{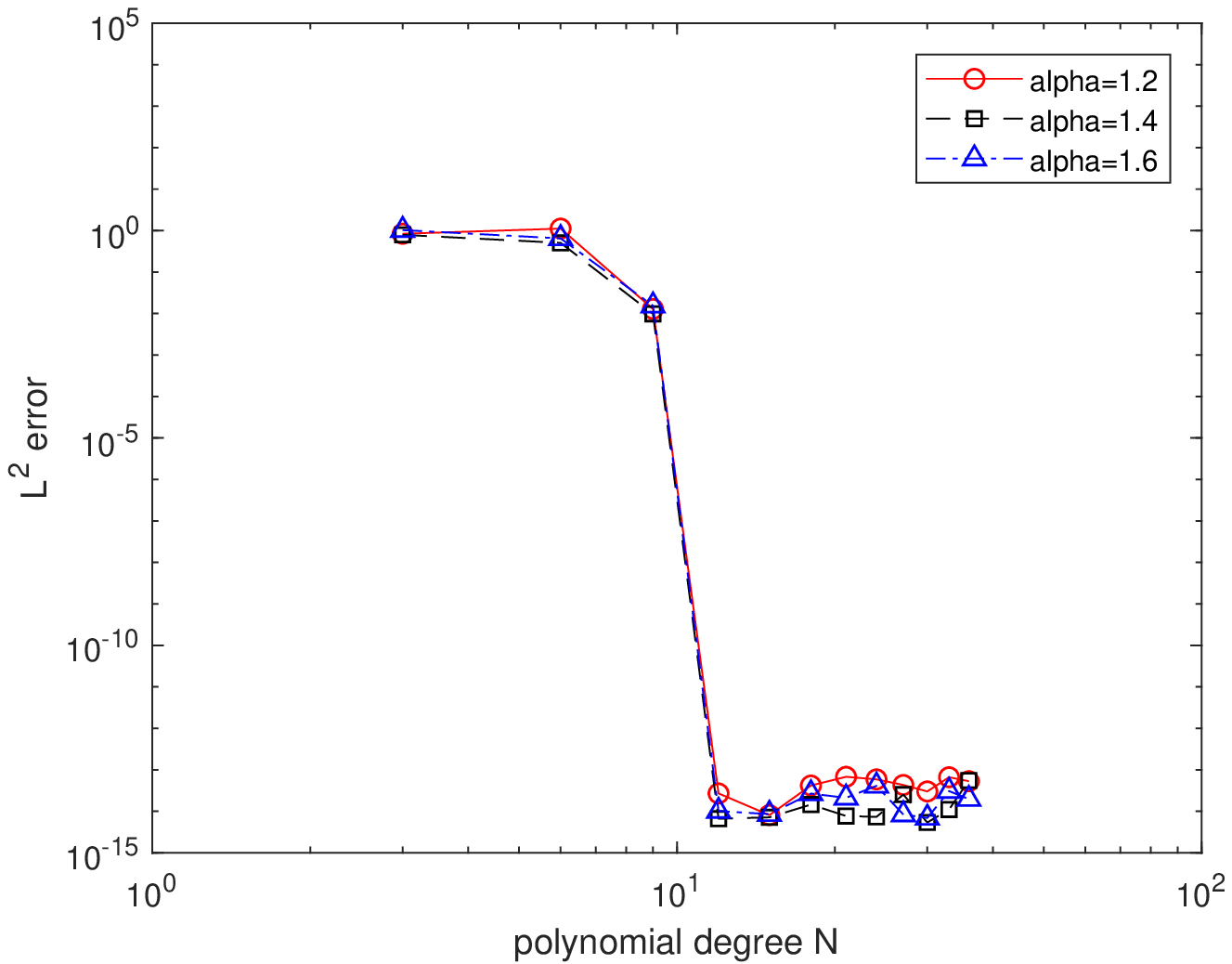}}
  \centerline{(b) Scheme 2}
\end{minipage}
\caption{The numerical $L_2$ errors Example \ref{exam3} vs the polynomial degree $N$ with $\alpha=1.2,1.4,1.6$ for Scheme 1 in (a) and Scheme 2 in (b).}
\label{fig3}
\end{figure}

\begin{example}\label{exam3}
  In this example, we verify the spectral convergence of Scheme 1 and Scheme 2 for the two-sided fractional diffusion equation with drift.
  Consider \eqref{main_equation} with $\bar{p}=\bar{q}=\frac{1}{2}$ and $\bar{d}=1$, i.e.,
  \begin{equation}\label{problem}
    -\frac{1}{2}({}_{-1}D_x^{\alpha,1}u(x) + {}_xD_1^{\alpha,1}u(x)) + u'(x)=h(x).
  \end{equation}
  For a given $\alpha$, the exact solutions for Scheme 1 and Scheme 2 are assumed to be with different forms, which are
\begin{equation}
\begin{array}{rl}
  &u_1(x)=(1+x)^{5+\frac{\alpha}{2}}(1-x)^5,\\[5pt]
  &u_2(x)=(1+x)^{5+\frac{\alpha-1}{2}}(1-x)^5,
\end{array}
\end{equation}
respectively.
\end{example}
 The numerical results for different $\alpha$ are shown in Figure \ref{fig3}, where spectral convergence can be observed when $N>10$ for both Scheme 1 and Scheme 2. We do not observe spectral convergence for Scheme 3. This might due to the three sub-equations with their individual trail functions we choose in Scheme 3. On the other hand, the advantage of Scheme 3 comes from its low condition number for all $\alpha\in(0,2)$, which will be illustrated in detail in the next example.

\begin{example}\label{exam7}
  Consider \eqref{main_equation} with $\bar{p}=\bar{q}=\frac{1}{2}$ and $\bar{d}=1$, i.e.,
  \begin{equation}\label{problem}
    -\frac{1}{2}({}_{-1}D_x^{\alpha}u(x) + {}_xD_1^{\alpha}u(x)) + u'(x)=h(x).
  \end{equation}
  We illustrate the condition number of the coefficient matrix for three schemes when solving problem \eqref{problem} with $\alpha=0.5,1.5,1.9$.
\end{example}

\begin{figure}[ht]
\begin{minipage}{0.31\linewidth}
  \centerline{\includegraphics[scale=0.3]{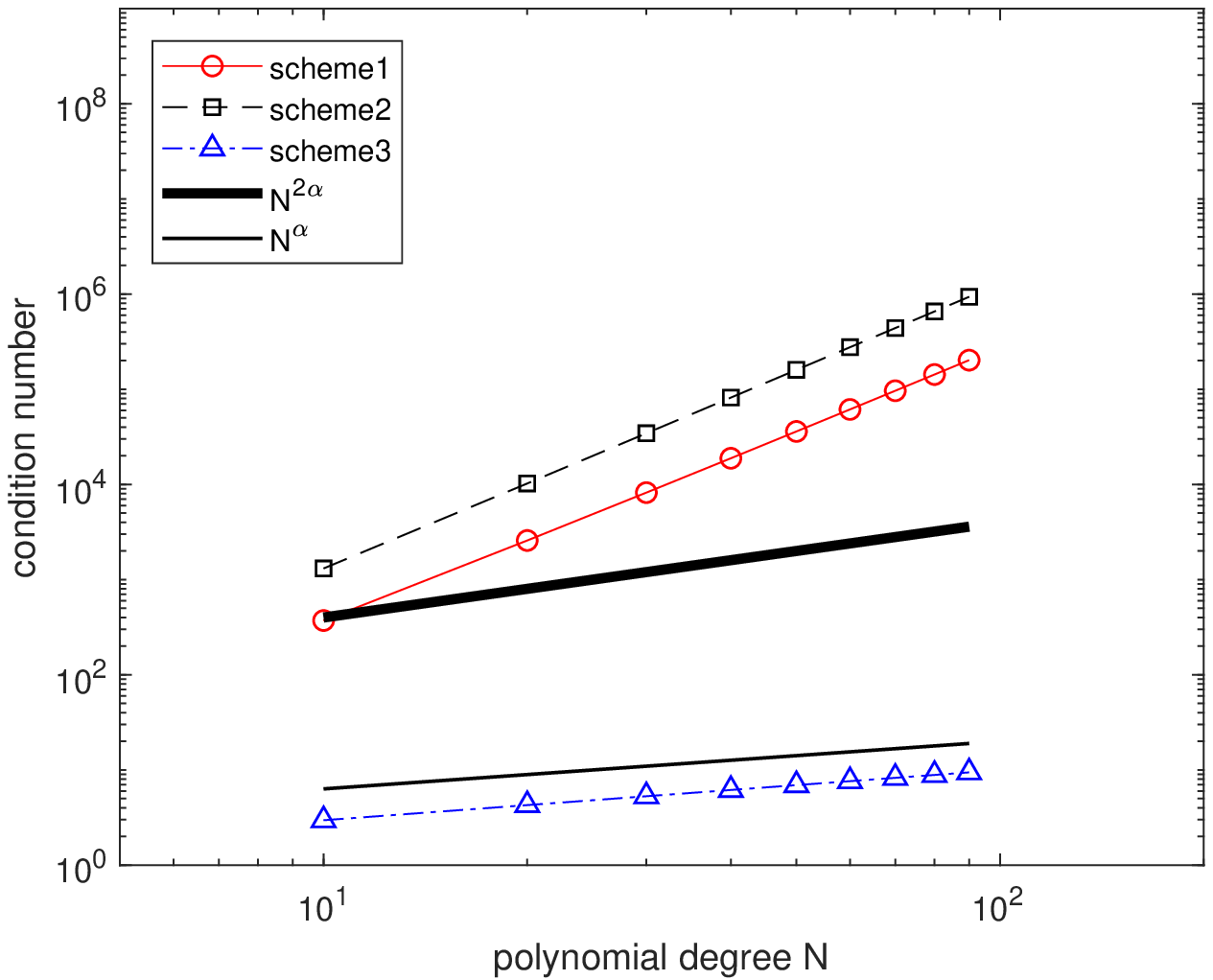}}
  \centerline{(a) $\alpha=0.5$}
\end{minipage}
\hfill
\begin{minipage}{0.31\linewidth}
  \centerline{\includegraphics[scale=0.3]{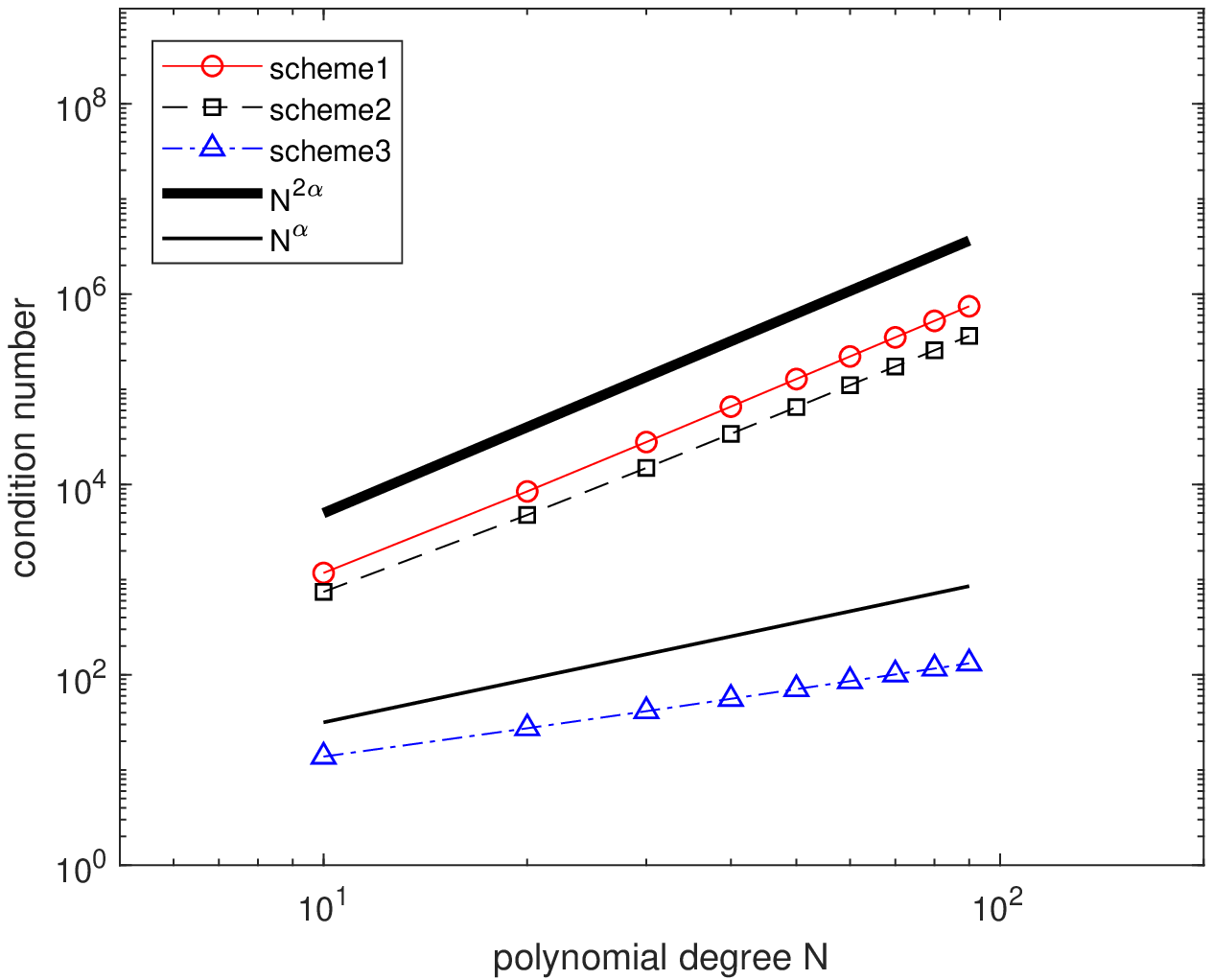}}
  \centerline{(b) $\alpha=1.5$}
\end{minipage}
\hfill
\begin{minipage}{0.31\linewidth}
  \centerline{\includegraphics[scale=0.3]{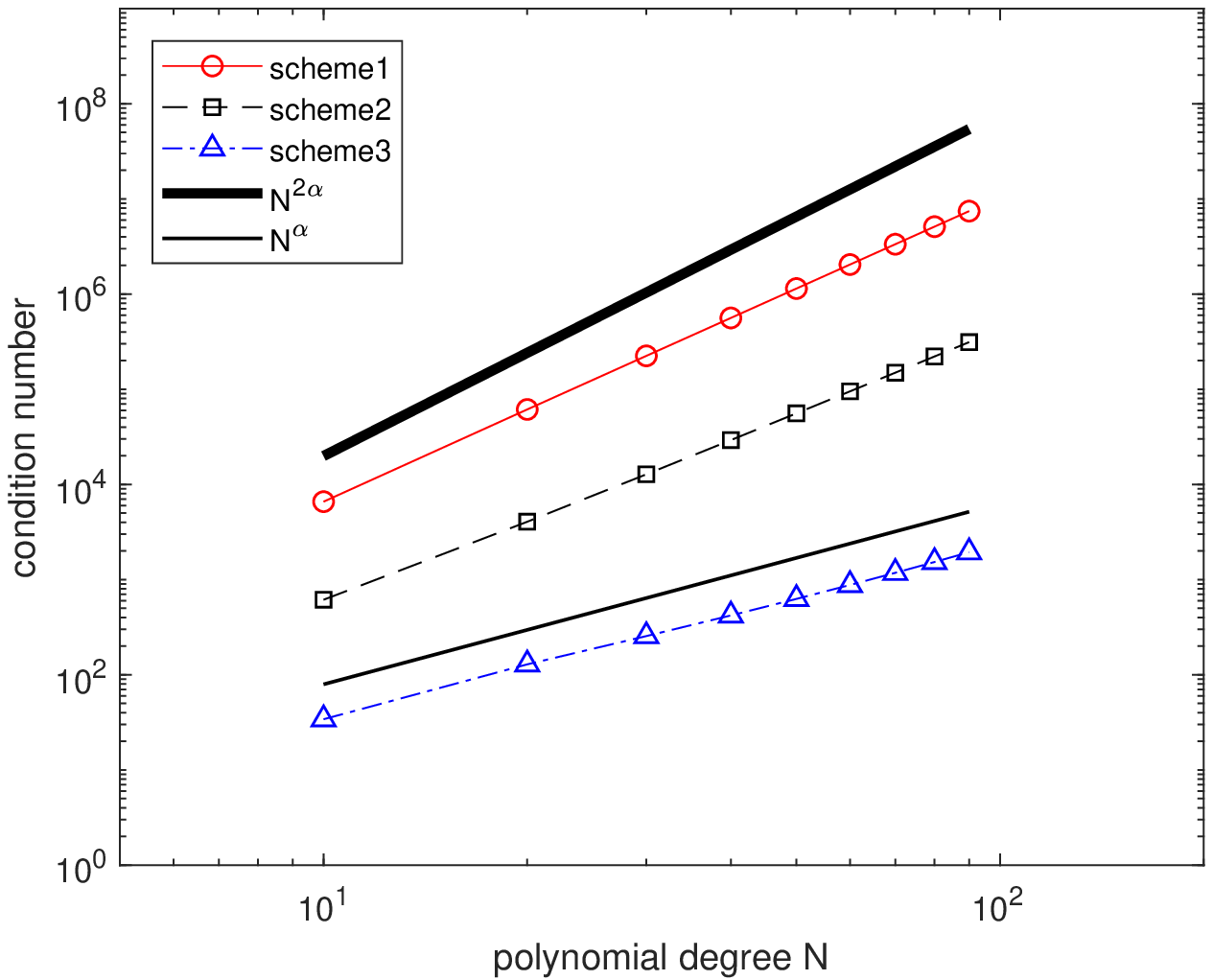}}
  \centerline{(c) $\alpha=1.9$}
\end{minipage}
\caption{The condition number versus polynomial degree $N$ of three schemes when solving Example \ref{exam7} with $\alpha=0.5,1.5,1.9$.}\label{fig7}
\end{figure}

As shown in Figure \ref{fig7}, the condition numbers of Scheme 1 and Scheme 2 grow as fast as $O(N^{2\alpha})$, and even faster when $\alpha=0.5$. But the condition numbers of Scheme 3 grow as $O(N^\alpha)$, which is much more moderate than the other two schemes.

\begin{figure}[ht]
\begin{minipage}{0.31\linewidth}
  \centerline{\includegraphics[scale=0.3]{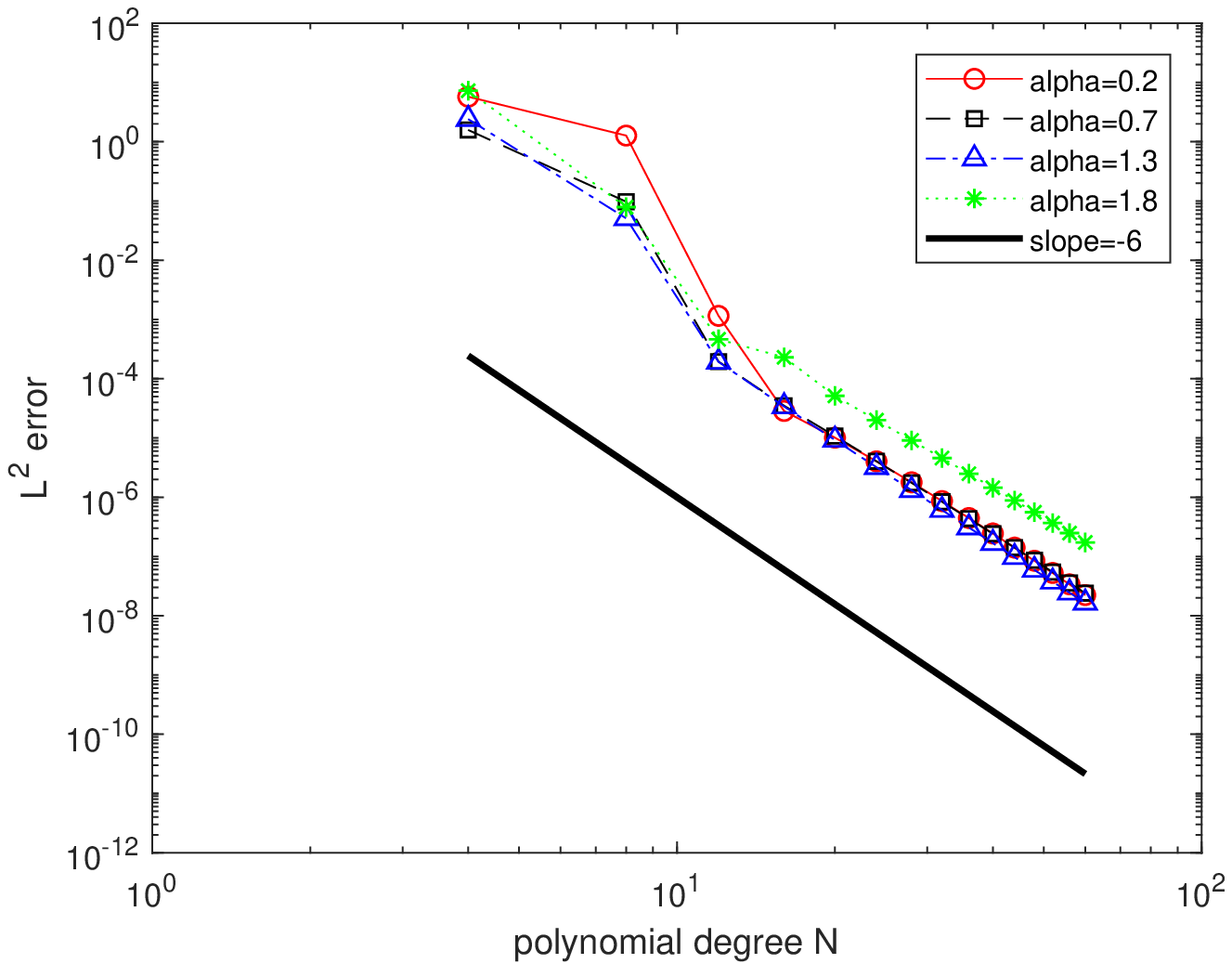}}
  \centerline{(a) $u_1=(1+x)^4(1-x)^3$}
\end{minipage}
\hfill
\begin{minipage}{0.31\linewidth}
  \centerline{\includegraphics[scale=0.3]{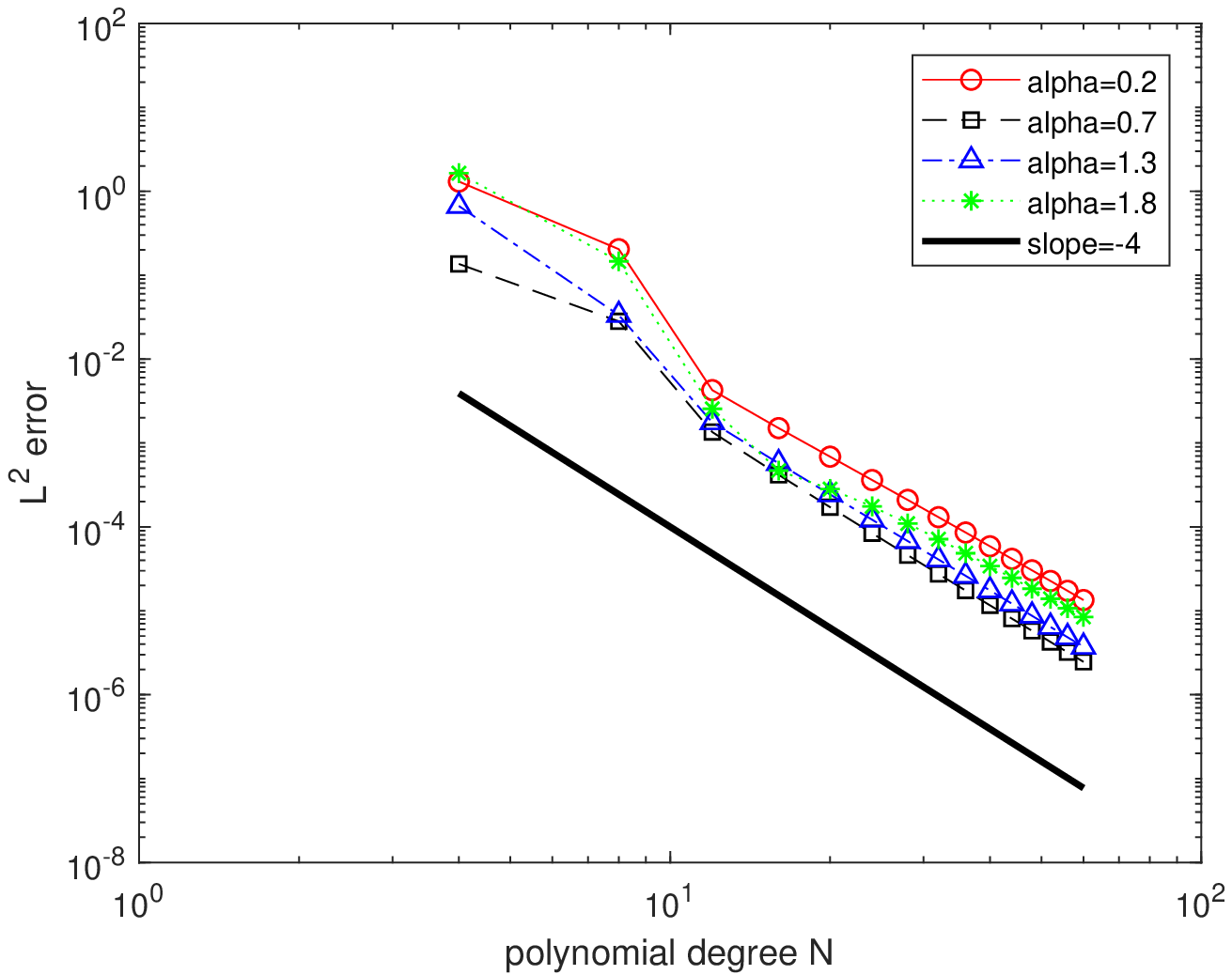}}
  \centerline{(b) $u_2=(1+x)^2(1-x)^4$}
\end{minipage}
\hfill
\begin{minipage}{0.31\linewidth}
  \centerline{\includegraphics[scale=0.3]{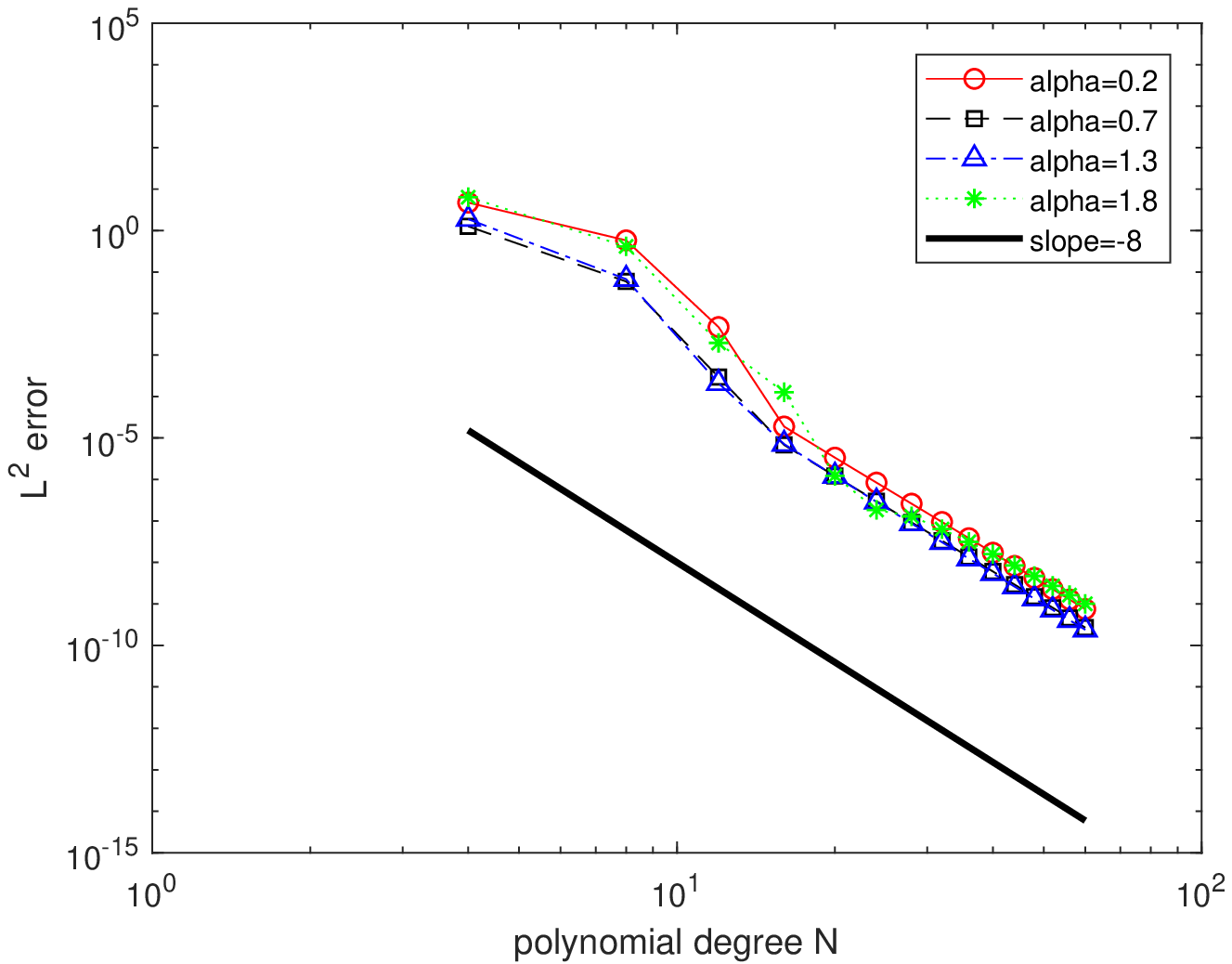}}
  \centerline{(c) $u_3=(1+x)^4(1-x)^4$}
\end{minipage}
\caption{The numerical $L_2$ errors of Example \ref{exam4} vs the polynomial degree N with Scheme 1.
(a), (b) and (c) show the errors decay algebraically when $\alpha=0.2,0.7,1.3,1.8$, with the exact solution being $u_1,u_2,u_3$ in \eqref{Scheme1ES}, respectively.}
\label{fig4}
\end{figure}

\begin{example}\label{exam4}
Besides the previous examples with special solutions, now we take some numerical tests with high regularities.
We use Scheme 1 to solve the same problem as Eq. \eqref{problem}. Three different exact solution with different regularities are assumed to be
\begin{equation}\label{Scheme1ES}
\begin{split}
  &u_1(x)=(1+x)^4(1-x)^3,\\
  &u_2(x)=(1+x)^2(1-x)^4,\\
  &u_3(x)=(1+x)^4(1-x)^4.
\end{split}
\end{equation}
\end{example}
The associated forcing term $h(x)$ cannot be given analytically. Instead, we compute $h(x)$ numerically at each Gauss quadrature nodes $x_i$.
We plot the $L_2$ error vs the polynomial degree $N$ for different values of $\alpha\in(0,2)$ in Figure \ref{fig4}, where $\alpha$ is taken to be $0.2,0.7,1.3,1.8$ and $N$ is from $4$ to $60$.
The errors show an algebraical decay and they are independent of the value of $\alpha$, only depend on the regularity of the exact solution $u$.

\begin{figure}[ht]
\begin{minipage}{0.31\linewidth}
  \centerline{\includegraphics[scale=0.3]{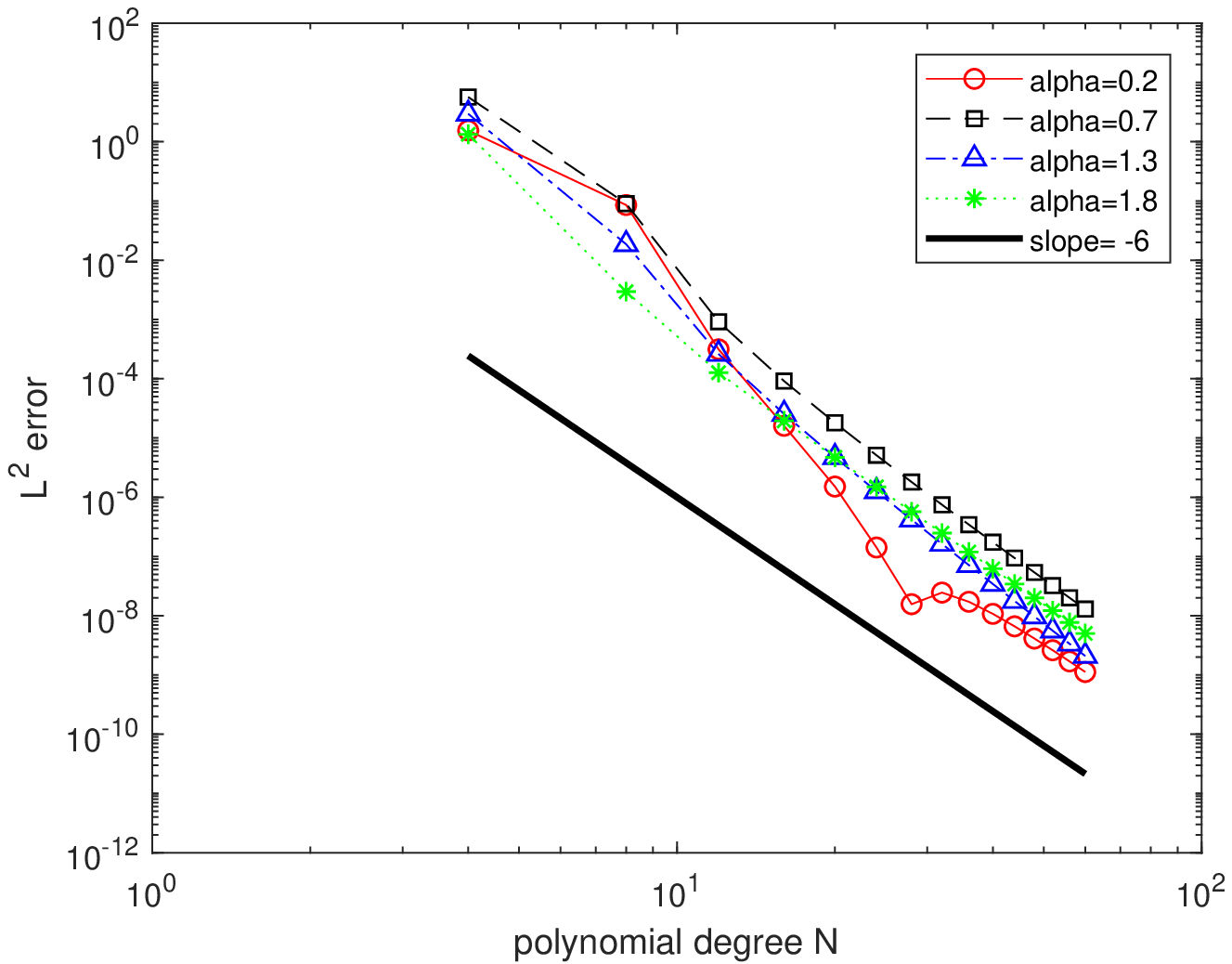}}
  \centerline{(a) $u_1=(1+x)^4(1-x)^3$}
\end{minipage}
\hfill
\begin{minipage}{0.31\linewidth}
  \centerline{\includegraphics[scale=0.3]{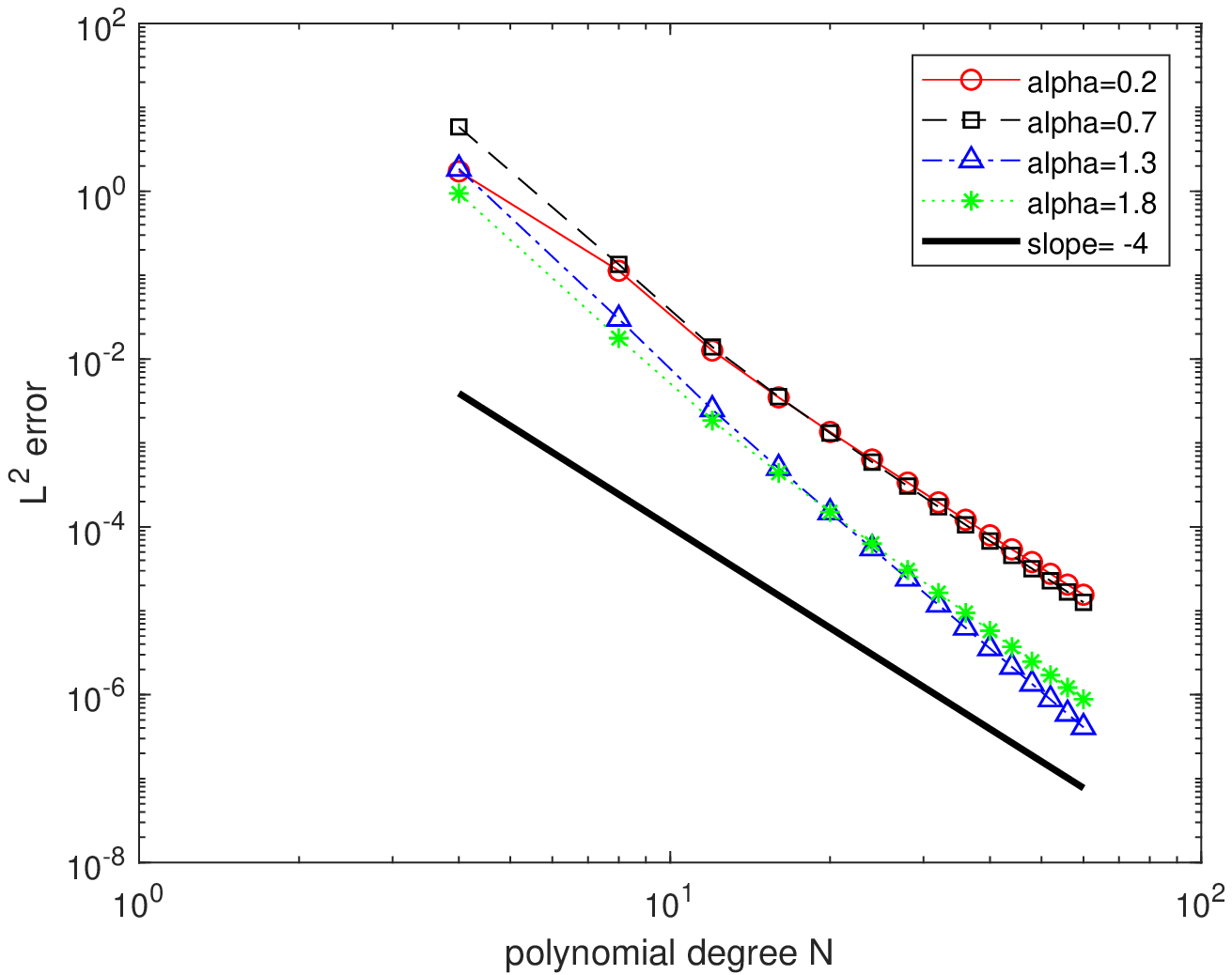}}
  \centerline{(b) $u_2=(1+x)^2(1-x)^4$}
\end{minipage}
\hfill
\begin{minipage}{0.31\linewidth}
  \centerline{\includegraphics[scale=0.3]{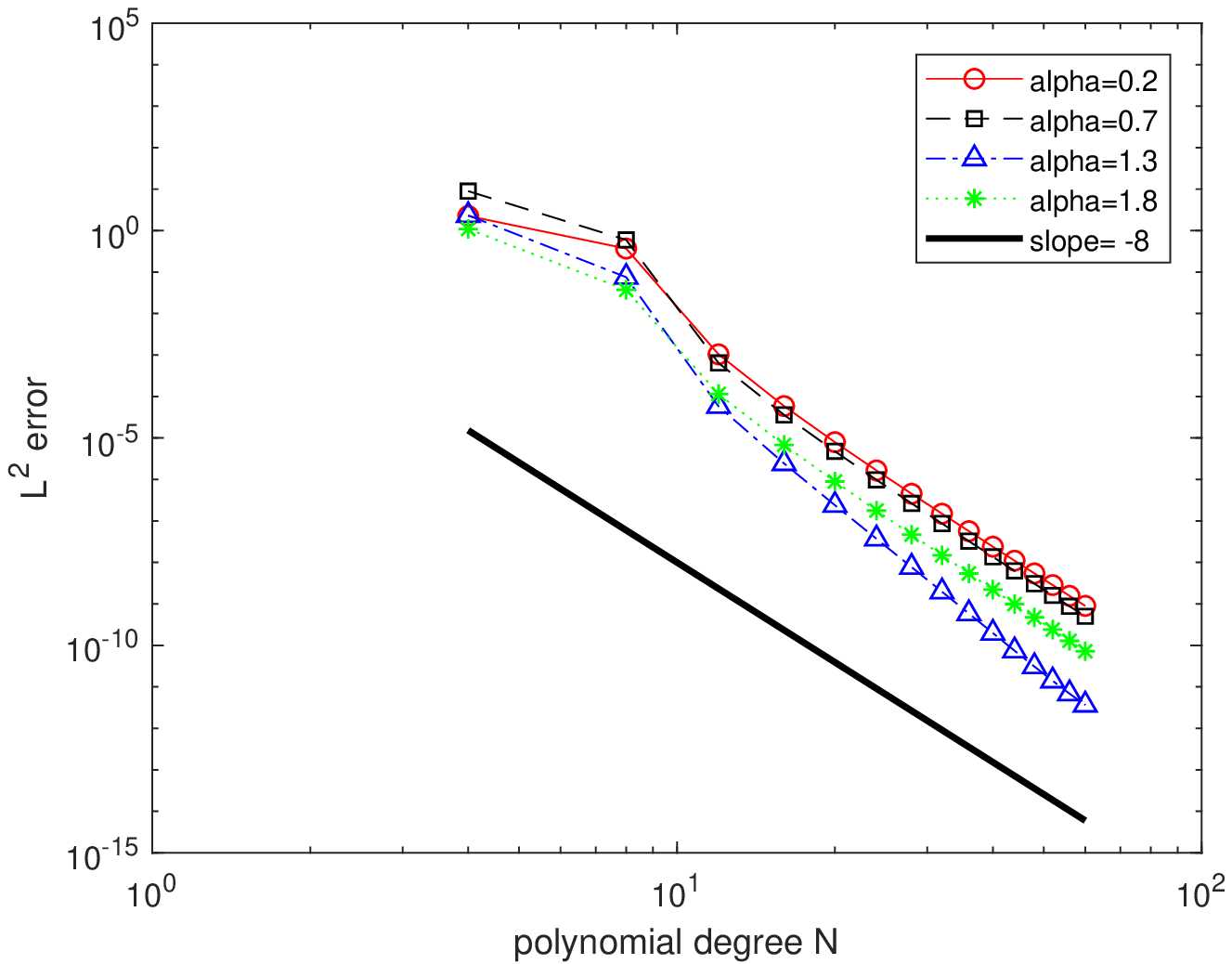}}
  \centerline{(c) $u_3=(1+x)^4(1-x)^4$}
\end{minipage}
\caption{The numerical $L_2$ errors of Example \ref{exam5} vs the polynomial degree N with Scheme 2.
(a), (b) and (c) show the errors decay algebraically when $\alpha=0.2,0.7,1.3,1.8$, with the exact solution being $u_1,u_2,u_3$ in \eqref{Scheme1ES}, respectively.}
\label{fig5}
\end{figure}

\begin{example}\label{exam5}
  We consider the same problem as Eq. \eqref{problem}, but using Scheme 2. The numerical results are shown in Figure \ref{fig5}.
\end{example}

\begin{figure}[ht]
\begin{minipage}{0.31\linewidth}
  \centerline{\includegraphics[scale=0.3]{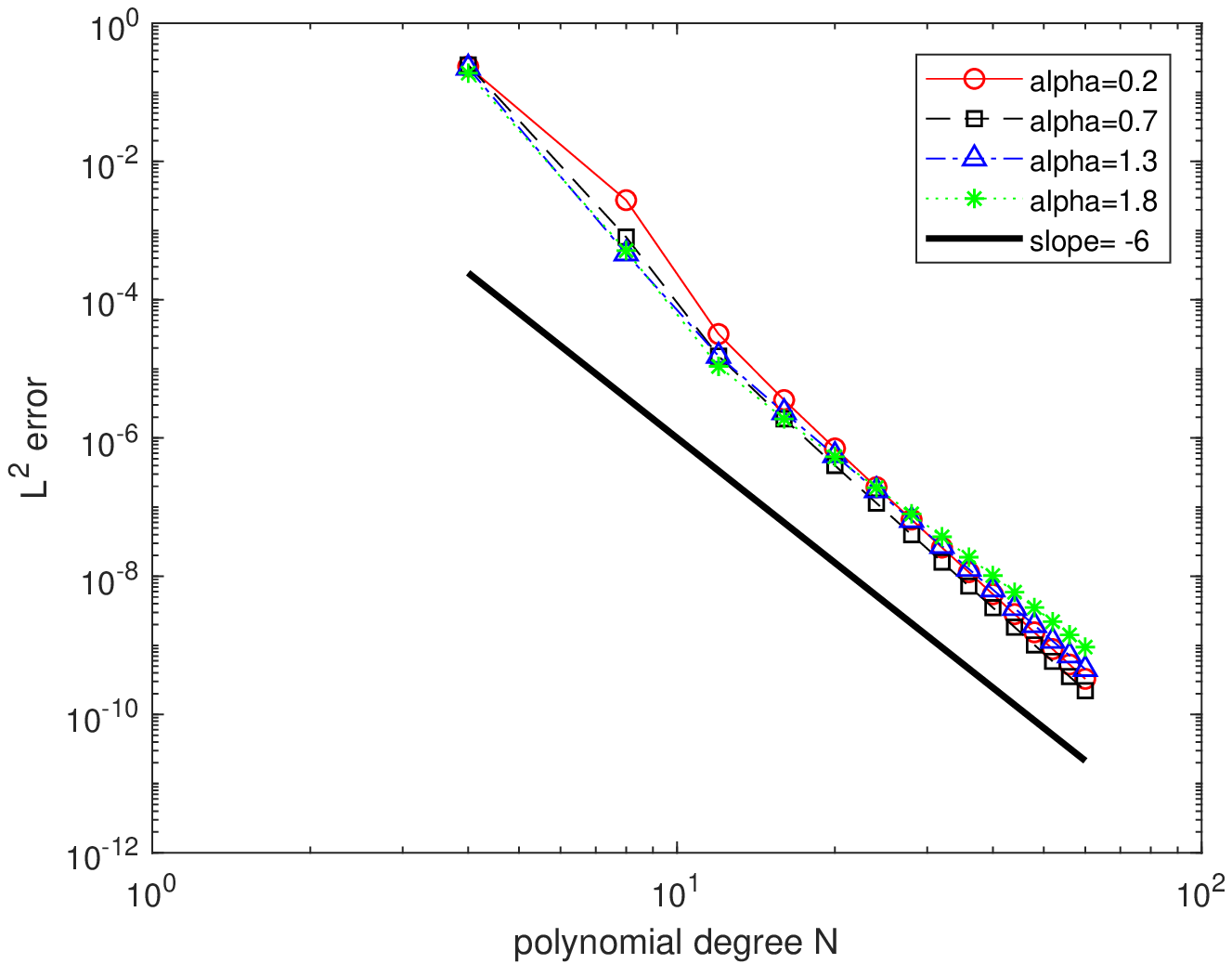}}
  \centerline{(a) $u_1=(1+x)^4(1-x)^3$}
\end{minipage}
\hfill
\begin{minipage}{0.31\linewidth}
  \centerline{\includegraphics[scale=0.3]{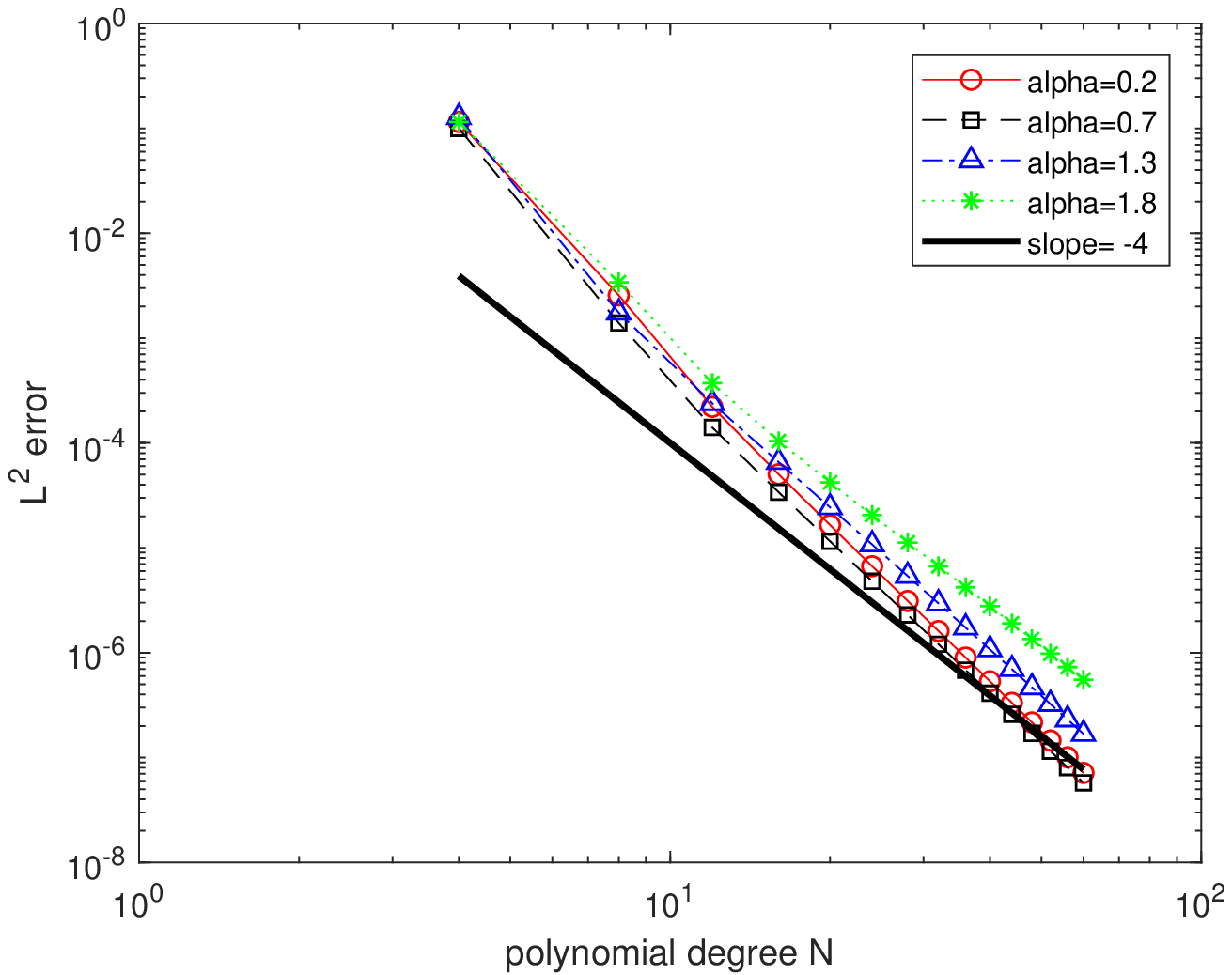}}
  \centerline{(b) $u_2=(1+x)^2(1-x)^4$}
\end{minipage}
\hfill
\begin{minipage}{0.31\linewidth}
  \centerline{\includegraphics[scale=0.3]{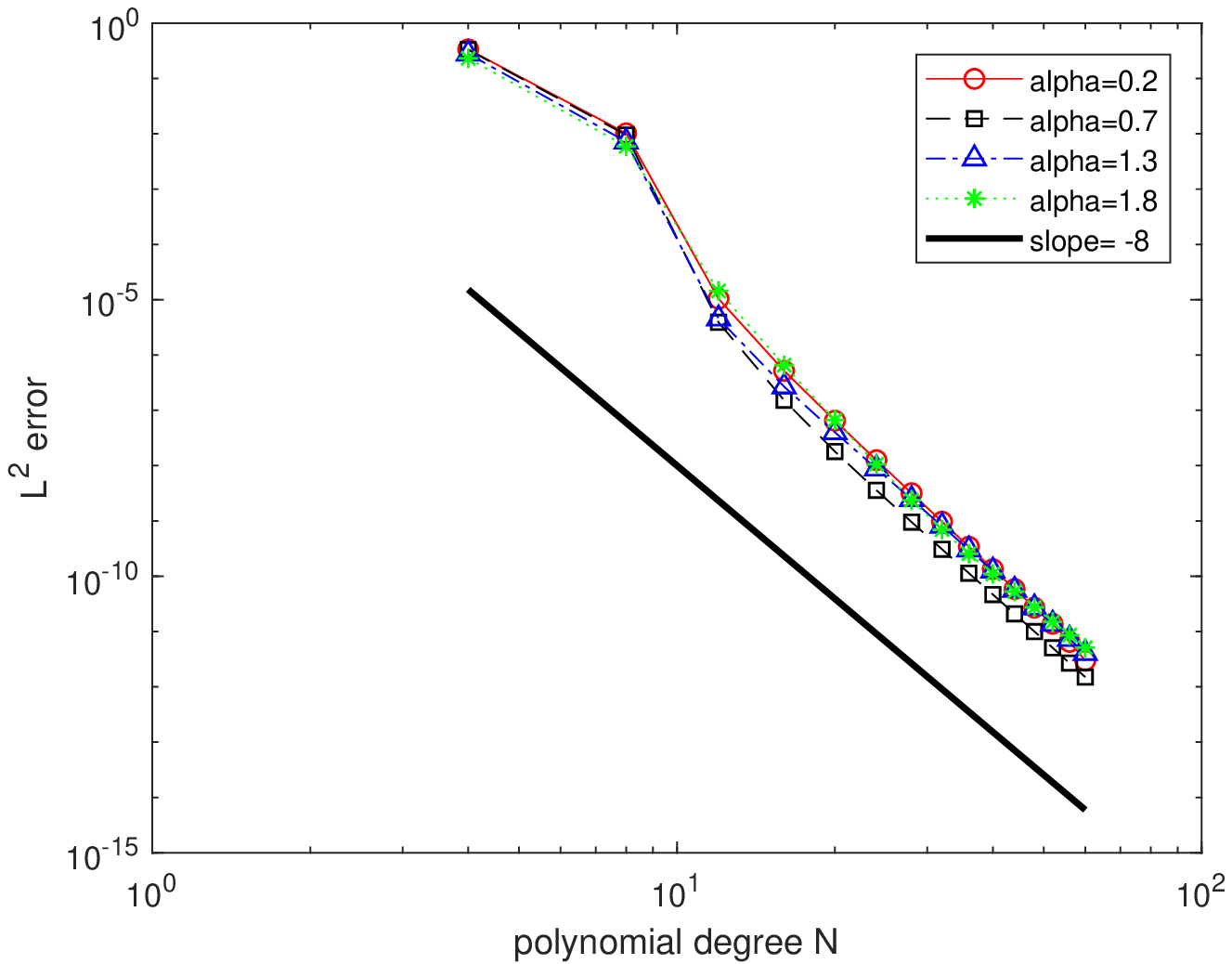}}
  \centerline{(c) $u_3=(1+x)^4(1-x)^4$}
\end{minipage}
\caption{The numerical $L_2$ errors of Example \ref{exam6} vs the polynomial degree N with Scheme 3.
(a), (b) and (c) show the errors decay algebraically when $\alpha=0.2,0.7,1.3,1.8$, with the exact solution being $u_1,u_2,u_3$ in \eqref{Scheme1ES}, respectively.}
\label{fig6}
\end{figure}

\begin{example}\label{exam6}
  We consider the same problem as Eq. \eqref{problem}, but using Scheme 3. The numerical results are shown in Figure \ref{fig6}.
\end{example}

Comparing the three Figures with three different schemes, we find that they are all effective to any $\alpha\in(0,2)$. For a general exact solution, their rates of convergence depend on the regularity of the solution on both sides.

\begin{example}\label{exam7}
Since the solution of the Eq. \eqref{main_equation} can be solved numerically with the three schemes we proposed in this paper, some interesting phenomena can be observed by analysing the numerical solution. One typical example is that if the source term $h(x)$ in Eq. \eqref{main_equation} taken to be $-1$, the solution $u(x)$ represents the mean first exit time of the particle starting at position $x\in(-1,1)$ when leaving the domain $[-1,1]$ \cite{Deng:17}.

Now, we firstly consider the mean first exit time of free diffusive particles at any given position $x\in(-1,1)$.
By taking $p=q=1/2$ and the source term $h(x)=\cos(\pi\alpha/2)$, we obtain the equivalent equation
\begin{equation}
  (-\Delta)^{\alpha/2}u(x)=-1.
\end{equation}
The numerical results for different $\alpha$ are shown in the left graph of Figure \ref{fig8}. One can observe that: the mean first exit time increases as $\alpha$ decreases; for each fixed $\alpha$, it costs more time for the particles in the middle part than those at near the boundary---all of these phenomenons are compatible with expectation. Next, we consider the effect of the drift term on the mean first exist time, by taking $d=\cos(\pi\alpha/2)$ which yields a drift to the left. Then the equation becomes
\begin{equation}
  (-\Delta)^{\alpha/2}u(x)- u'(x)=-1.
\end{equation}
The corresponding results are demonstrated in the right graph of Figure \ref{fig8}. In this case, under the effect of the drift, the particles are more likely to leave the domain from the left side. In other words, it takes more time for the particles at the right part to leave the domain. If one particle starts at the very right part (near the right boundary), however, the diffusion behavior works and makes the particle leave the right boundary in a moment time.

\begin{figure}[ht]
\begin{minipage}{0.45\linewidth}
  \centerline{\includegraphics[scale=0.4]{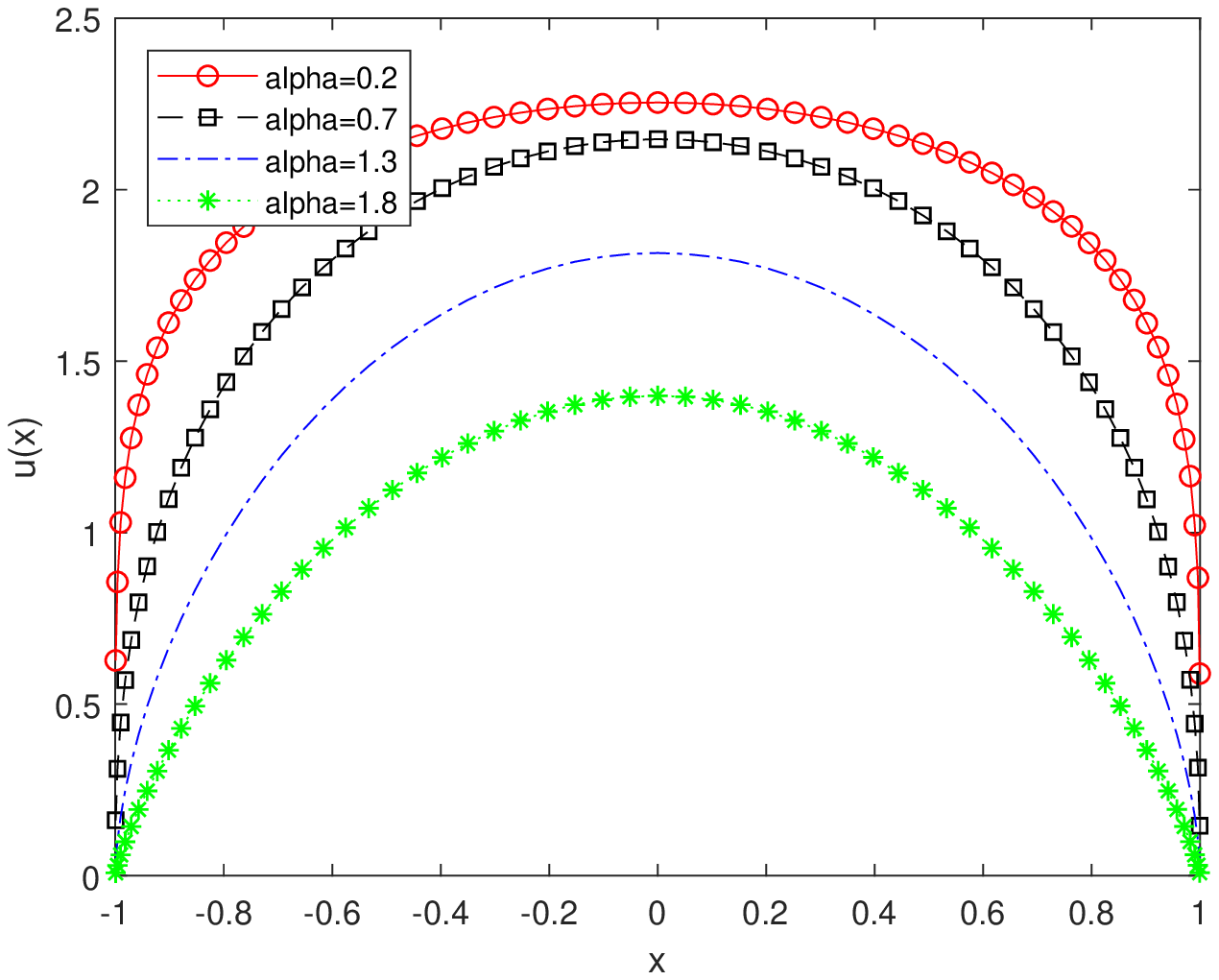}}
  \centerline{(a) $d=0$}
\end{minipage}
\hfill
\begin{minipage}{0.45\linewidth}
  \centerline{\includegraphics[scale=0.4]{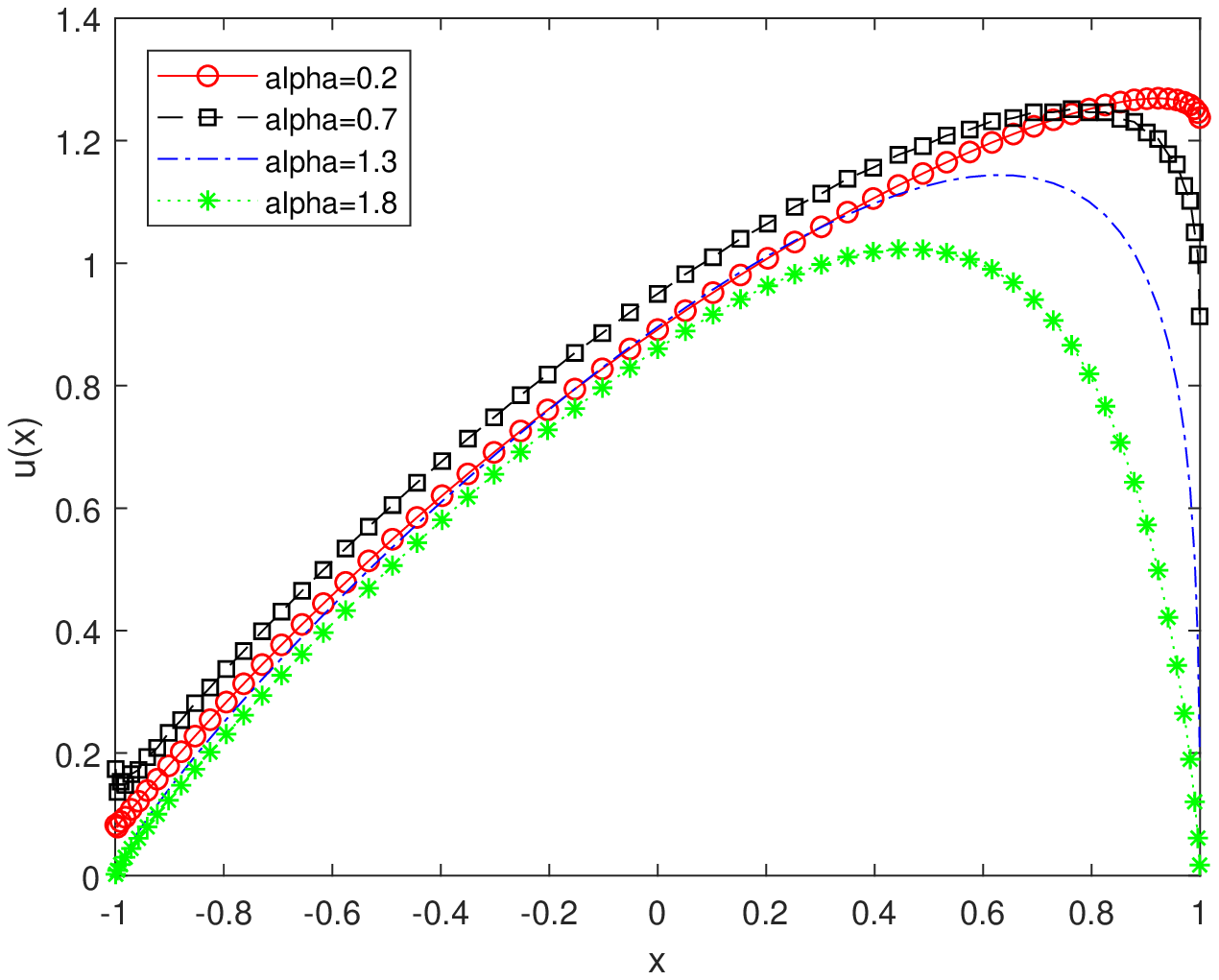}}
  \centerline{(b) $d=\cos(\pi\alpha/2)$}
\end{minipage}
\caption{The graph of solution $u(x)$ of Example \ref{exam7} with $\alpha=0.2,0.7,1.3,1.8$. Here, $u(x)$ represents the mean first exit time of the particle starting at position $x\in(-1,1)$ when leaving this domain.}
\label{fig8}
\end{figure}

\end{example}

\section{Conclusion}\label{section:5}
In this paper, we discuss spectral approximations in the weak sense for solving a two-sided fractional differential equation with drift, in which the fractional operators are physically well-defined \cite{Ervin:16}.
Three kinds of spectral formulae, namely Galerkin spectral formulation, Petrov-Galerkin spectral formulation, and mixed Galerkin spectral formulation, are proposed step by step. Then their corresponding spectral Galerkin schemes are derived. The significant advantage of the mixed Galerkin spectral scheme is that its condition number grows as $O(N^\alpha)$, compared with the other two schemes, whose condition numbers grow as fast as $O(N^{2\alpha})$.
We compare these three kinds of schemes through several numerical experiments. All of them turn out to be effective for different problems, especially also for the fractional Laplacian with generalized Dirichlet boundary conditions, the fractional order of which is $\alpha\in(0,2)$, not only having to be limited in $(1,2)$. What is more, considering the physical meanings of the fractional differential equation with drift, one interesting physical quantity, mean first exit time, is computed and discussed in this paper.
More related theoretical analysis will be discussed in our future work.
\bibliographystyle{amsplain}

\end{document}